\documentclass[11pt]{article}
\usepackage{amscd}      
\usepackage{amssymb}
\usepackage[intlimits]{amsmath}
\usepackage{xypic}      
\LaTeXdiagrams          
\usepackage[all,v2]{xy}
\xyoption{2cell}
\UseAllTwocells
\xyoption{frame}
\CompileMatrices

\usepackage[all]{xy}

\usepackage{theorem}

\newtheorem{prop}{Proposition}[section]

\newtheorem{lem}[prop]{Lemma}

\newtheorem{cor}[prop]{Corollary}

\newtheorem{them}[prop]{Theorem}

\theorembodyfont{\upshape}

\newtheorem{defn}[prop]{Definition}

\newtheorem{rmk}[prop]{Remark}

\newtheorem{numex}[prop]{Example}

\newtheorem{numrk}[prop]{Remark}

\newenvironment{pf}{\begin{trivlist}\item[]{\sc Proof.}}%
            {\nolinebreak $\Box$ \end{trivlist}}

            {\nolinebreak $\Box$ \end{trivlist}}

\newcommand{\noprint}[1]{}

\newcommand{\equal}{=}

\renewcommand{\tilde}{\widetilde}

\newcommand{\toto}{\rightrightarrows}

\newcommand{\com}{{\scriptscriptstyle\bullet}}

\newcommand{\Gg}{{\mathfrak g}}
\newcommand{\Hh}{{\mathfrak h}}

\newcommand{\nn}{{\mathbb N}}

\newcommand{\rr}{{\mathbb R}}

\newcommand{\del}{\partial}

\newcommand{\pr}{\mathop{\rm pr}\nolimits}

\newcommand{\id}{\mathop{\rm id}\nolimits}

\newcommand{\smalcirc}{\mbox{\tiny{$\circ $}}}

\newcommand{\ldiag}[1]%
       {\makebox[0cm]{${\scriptstyle#1}\downarrow\phantom{\scriptstyle#1}$}}

\newcommand{\ldiagup}[1]%
       {\makebox[0cm]{${\scriptstyle#1}\uparrow\phantom{\scriptstyle#1}$}}

\newcommand{\rdiag}[1]%
       {\makebox[0cm]{$\phantom{\scriptstyle#1}\downarrow{\scriptstyle#1}$}}

\newcommand{\sediagr}[1]%
       {\makebox[0cm]{$\phantom{\scriptstyle#1}\searrow{\scriptstyle#1}$}}

\newcommand{\nediagr}[1]%
       {\makebox[0cm]{$\phantom{\scriptstyle#1}\nearrow{\scriptstyle#1}$}}

\newcommand{\rdiagup}[1]%
       {\makebox[0cm]{$\phantom{\scriptstyle#1}\uparrow{\scriptstyle#1}$}}

\newcommand{\swdiag}[1]%
       {\makebox[0cm]{$\phantom{\scriptstyle#1}\swarrow{\scriptstyle#1}$}}

\newcommand{\sediag}[1]%
       {\makebox[0cm]{${\scriptstyle#1}\searrow\phantom{\scriptstyle#1}$}}

\newcommand{\nediag}[1]%
       {\makebox[0cm]{${\scriptstyle#1}\nearrow\phantom{\scriptstyle#1}$}}

\newcommand{\iso}{\stackrel{\sim}{\rightarrow}}

\newcommand{\doublearrowstack}[2]%
 {{{{\scriptstyle#1}\atop{\textstyle\longrightarrow}}\atop{{\textstyle\longright
arrow}\atop{\scriptstyle#2}}}}
\newcommand{\rightleftarrowstack}[2]%
 {{{{\scriptstyle#1}\atop{\textstyle\longrightarrow}}\atop{{\textstyle\longlefta
rrow}\atop{\scriptstyle#2}}}}
\newcommand{\leftrightarrowstack}[2]%
 {{{{\scriptstyle#1}\atop{\textstyle\longleftarrow}}\atop{{\textstyle\longrighta
rrow}\atop{\scriptstyle#2}}}}

\newcommand{\comdia}[9]{%
\begin{array}{ccc}
#1 & \stackrel{#2}{\longrightarrow} & #3 \\
\ldiag{#4} & #5 & \rdiag{#6} \\
#7 & \stackrel{#8}{\longrightarrow} & #9
\end{array}}

\newcommand{\comtri}[6]{%
\begin{array}{ccc}
#1 & \stackrel{#2}{\longrightarrow} & #3         \\
   & \sediag{#4}                    & \rdiag{#5} \\
   &                                & #6
\end{array}}

\newcommand{\comtriup}[6]{%
\begin{array}{ccc}
#1 & \stackrel{#2}{\longrightarrow} & #3         \\
   & \sediag{#4}                    & \rdiagup{#5} \\
   &                                & #6
\end{array}}

\newcommand{\overtoparrow}%
{\makebox[0cm]{\beginpicture
\setcoordinatesystem units <.8cm,.4cm> point at 0 0
\setplotarea x from -3 to 3, y from 0 to 1
\setquadratic
\plot -3 0 0 1 3 0 /
\put{\vector(3,-1){0}}[Bl] at 3 0
\endpicture}}

\newcommand{\underbottomarrow}%
{\makebox[0cm]{\beginpicture
\setcoordinatesystem units <.8cm,.4cm> point at 0 0
\setplotarea x from -3 to 3, y from 0 to 1
\setquadratic
\plot -3 1 0 0 3 1 /
\put{\vector(3,1){0}}[Bl] at 3 1
\endpicture}}

\newcommand{\ses}[5]%
{0\longrightarrow#1\stackrel{#2}{ \longrightarrow}#3\stackrel{#4}{
\longrightarrow}#5\longrightarrow0}

\newcommand{\dt}[6]%
{#1\stackrel{#2}{\longrightarrow}#3 \stackrel{#4}{\longrightarrow}#5
\stackrel{#6}{\longrightarrow} #1[1]}

\newcommand{\cat}[1]%
{(\mbox{\rm #1})}

\newcommand{\gm}{\Gamma}

\newcommand{\lon }{\longrightarrow }

\newcommand{\be }{\begin{eqnarray*}}

\newcommand{\ee }{\end{eqnarray*}}

\newcommand{\source}{source }
\newcommand{\target}{target }
\newcommand{\ttheta}{\eta}
\newcommand{\so}{\beta}
\newcommand{\ta}{\alpha}
\newcommand{\per }{\backl }
\newcommand{\backl}{\mathbin{\vrule width1.5ex height.4pt\vrule height1.5ex}}
 \newcommand{\ii}{i}
\newcommand{\half}{\frac{1}{2}}

\newcommand{\Diff}{\mbox{Diff}}
\newcommand{\Ker}{ker}
\newcommand{\la}{\langle}
\newcommand{\ra}{\rangle}
\newcommand{\degeneracy}{face }
\newcommand{\face}{degeneracy }
\newcommand{\tsigma}{\tilde{\tau}}
\newcommand{\ttau}{\tilde{\sigma}}
\newcommand{\cw}{w}
\newcommand{\eepsilon}{\eta}
\newcommand{\eeta}{\epsilon}
\def\gpd{\,\lower1pt\hbox{$\longrightarrow$}\hskip-.24in\raise2pt
             \hbox{$\longrightarrow$}\,}

\newcommand{\LL}{{\mathcal L}}

\let\Tilde=\widetilde

\let\Vec=\overrightarrow

\def\chigh{{\raise1.5pt\hbox{$\chi$}}}
\let\phi=\varphi
\def\til0{\Tilde{0}}

\title{Chern-Weil map for principal bundles over groupoids}

\author{ Camille Laurent-Gengoux \\
 D\'epartement de Math\'ematiques \\ 
 Universit\'e de Poitiers \\
SP2MI, Boulevard Marie et Pierre Curie, \\
86962 Futuroscope-Chasseneuil Cedex, France. \\
{\sf email: Camille.Laurent@math.univ-poitiers.fr }
\\\\
 Jean-Louis Tu\\
Laboratoire de Math\'ematiques et Applications de Metz\\
Universit\'e de Metz\\
 ISGMP, B\^atiment A, Ile du Saulcy\\
 57000 Metz, France\\   
 {\sf email:  tu@univ-metz.fr}
 \\\\
AND\\
Ping Xu \thanks{Research partially supported by NSF
       grant DMS03-06665. }\\
        Department of Mathematics\\
         Pennsylvania State University \\
         University Park, PA 16802, USA\\
{\sf email: ping@math.psu.edu }
}

\begin{document}

\sloppy

\maketitle

\begin{abstract}
The theory of principal $G$-bundles over a Lie groupoid is an important one
unifying  various types of principal $G$-bundles, including those over
manifolds, those over orbifolds, as well as equivariant principal $G$-bundles.
In this paper, we study differential geometry of these objects, including
connections and holonomy maps. We also introduce a Chern-Weil map for these
principal bundles and prove that the characteristic classes obtained coincide
with the universal characteristic classes.

As an application, we recover the equivariant Chern-Weil map of Bott-Tu. We also
obtain an explicit chain map between the Weil model and the simplicial model of
equivariant cohomology which reduces to the Bott-Shulman map $S(\Gg^*)^G \to
H^*(BG)$ when the manifold is a point.
\end{abstract}

{\small \tableofcontents}

\section{Introduction}

A remarkable, and well-known, theory based on principal $G$-bundles over a
manifold is the theory of Chern-Weil \cite{Chern, Weil} which constructs
characteristic classes geometrically in terms of differential forms. Recently,
there has been an increasing  interest in other types
 of principal bundles, for
instance, principal bundles over orbifolds or equivariant principal bundles. It
is therefore important to develop a theory of characteristic classes for these
more general settings.
For equivariant principal bundles, this has been done by many authors
\cite{BGV, BT, Jeffrey}, who
introduced a Chern-Weil map  in which the de Rham cohomology
 groups in the usual Chern-Weil theory were replaced by
equivariant cohomology (using the Weil model or  the Cartan model).
More recently, Alekseev-Meinrenken gave a construction of Chern-Weil map for
non-commutative $\Gg$-differential algebras \cite{AM03}.

In this paper, we develop Chern-Weil theory for principal $G$-bundles over Lie
groupoids. The notion of a groupoid is a standard generalization of the concepts
of \emph{spaces} and \emph{groups}, in which spaces and groups are treated on an
equal footing. Simplifying somewhat, one could say that a groupoid is a mixture
of a space and a group, with properties of both interacting in a delicate way.
In a certain sense, groupoids provide a uniform framework for many different
geometric objects. For instance, when a Lie group $G$ acts on a manifold $M$ the
transformation groupoid $G\times M \toto M$ can be used to replace the quotient
space. On the other hand, an orbifold can be represented by a proper \'etale
groupoid \cite{Moe}. Thus the notion of principal bundles over groupoids unifies
various existing notions of principal bundles.

A Lie groupoid $\gm \toto \gm_0$ gives rise to a simplicial manifold $\gm_\com$
in a natural fashion. Thus one can use the double complex $\Omega(\gm_\com)$
associated to the simplicial manifold $\gm_\com$ as a model for differential
forms. Its cohomology $H^*_{dR}(\gm_\com)$ is isomorphic to the cohomology of
the classifying space $B\gm$ \cite{Segal}. 
It is well known that $H^*_{dR}(\gm_\com)$ becomes a ring
under the cup product \cite{Dupont}.
 For a Lie group, this is the model studied by
Bott-Shulman in connection with characteristic classes \cite{Bott,BSS,BS}. Under this
model, for a principal $G$-bundle $P \to \gm_0$ over a Lie groupoid $\gm \toto
\gm_0$, the relevant notion involved is a pseudo-connection, {\it i.e.} a
connection $1$-form $\theta \in \Omega^1(P) \otimes \Gg$ for the $G$-bundle $P
\to \gm_0$ (forgetting the groupoid structure). The total pseudo-curvature is a
degree $2$element in $\Omega(Q_\com) \otimes \Gg$ consisting of two parts. The
first part is the usual curvature $2$-form $\Omega \equal d\theta +
\half[\theta,\theta] \in \Omega^2(P) \otimes \Gg$, while the second part
$\partial\theta \in \Omega^1(Q) \otimes \Gg$ measures the failure of $\theta$ to
be $\gm$-basic (Proposition \ref{prop:existeuneconnection}). 
Here $Q \toto P$ denotes the transformation groupoid associated
 to the $\gm$-action on $P$, {\it i.e.} $Q \equal \gm \times_{\gm_0} P$,
and $\partial :\Omega^1 (P) \otimes \Gg \to \Omega^1 (Q) \otimes \Gg$
is the boundary map  with respect to this transformation groupoid.
 Comparing this with the usual picture,
it is useful to note that the differential stack  corresponding to the
groupoid $Q\toto P$ is indeed a principal $G$-bundle over the differential stack of $\gm$. We prove

\vskip.2in

\noindent{\bf Theorem A.} {\it
\begin{enumerate}
\item Associated to any pseudo-connection $\theta \in \Omega^1(P) \otimes \Gg$,
on the cochain level there is a canonical map
    \[ z_\theta: S(\Gg^*)^G \to Z^*_{dR}(\gm_\com), \]
called the Chern-Weil map, where $Z^*_{dR}(\gm_\com)$ is the space of closed
forms. On the level of cohomology, $z_\theta$ induces an algebra homomorphism
    \[ w_\theta: S(\Gg^*)^G \to H^*_{dR}(\Gamma_\com), \]
which is independent of the choice of the pseudo-connection and thus may be
denoted by $w_P$. Moreover, $z_\theta$ is completely determined by the total
pseudo-curvature (Proposition \ref{pro:5.11}).
\item If $\phi$ is a strict homomorphism from $\gm' \toto \gm_0'$ to $\gm \toto
\gm_0$, then the following diagrams
    \[ \begin{array}{ccccccc}
        S({\Gg}^*)^G & \stackrel{z_{\phi^* \theta}}{\to} &
            Z_{dR}^*(\gm_{\com}') & & S({\Gg}^*)^G & \stackrel{w_{P'}}{\to} &
            H^*_{dR}(\gm_{\com}') \\
        & \stackrel{\searrow}{z_{\theta}}& \uparrow \phi^* &\mbox{and} & &
            \stackrel{\searrow}{w_P} & \uparrow \phi^* \\
        & & Z_{dR}^*(\gm_{\com}) & & & & H^*_{dR}(\gm_{\com})
    \end{array} \]
commute, where $P' \equal \phi^* P$ is the pull-back of $P$ via $\phi$.
\end{enumerate}}

\vskip.2in

As an application, we show that  this construction reduces to various existing
constructions in the literature in specific cases.

\vskip.2in

\noindent{\bf Theorem B.} {\it
\begin{enumerate}
\item If $\Gamma$ is a manifold $M$, then $w_P: S(\Gg^*)^G \to H^*_{dR}(M)$
reduces to the usual Chern-Weil map \cite{guillemin}.
\item If $\Gamma$ is a Lie group $G$, $P$ the $G$-bundle $G \to \cdot $, and
$\theta$ the left Maurer-Cartan form, then $z_\theta: S(\Gg^*)^G \to
Z_{dR}^*(G_\com)$ and $w_P: S(\Gg^*)^G \to H^*_{dR}(G_\com)$ coincide with the
Bott-Shulman maps \cite{Bott,BSS}.
\end{enumerate}}

\vskip.2in

Another interesting application is the case of equivariant principal
$G$-bundles. Let $P \to M$ be an $H$-equivariant principal $G$-bundle. In
\cite{BT}, Bott-Tu introduced a Chern-Weil map, with values in the Weil model,
associated to any $H$-invariant connection $\theta$ on $P \to M$. More
precisely, they constructed an $H$-basic connection $\Xi$ on the
$G$-differential algebra $W(\Hh) \otimes \Omega(P)$, which induces a map
$z_{BT}: S(\Gg^*)^G \to \big(W(\Hh) \otimes \Omega(M)\big)^{H-basic}$. On the
other hand, $P$ can be considered as a $G$-bundle over the transformation
groupoid $H \times M \toto M$, and $\theta \in \Omega(P)\otimes \Gg$ as a
pseudo-connection. Our construction in Section \ref{sec:CWgroupoides} induces a
map $z_{\theta}:S(\Gg^* )^G \to Z((H\times M)_\com)$. We have

\vskip.2in

\noindent{\bf Theorem C.} {\it
The following diagram  commutes
    \[ \comtri{S(\Gg^* )^G}{z_{BT}}
        {Z^*\big(W(\Hh) \otimes \Omega(M)\big)^{H-basic}}{z_\theta}{K}
        {Z((H \times M)_\com)} \]
where $K: W(\Hh) \otimes \Omega(M) \to \Omega((H\times M)_\com)$ is the natural
chain map between the Weil model and the simplicial model as described in
Proposition \ref{pro:5.17}.}

\vskip.2in

There is a conceptually simpler way to think of characteristic classes, namely
via  universal characteristic classes. A principal $G$-bundle $P$ over $M$
induces a map from $M$ to $BG$ (the classifying space of $G$), which in turn
induces a homomorphism of cohomology groups $H^*(BG) \to H^*(M)$. Composing the
Bott-Shulman map with this map, one gets a map $S(\Gg^*)^G \to H^*(M)$. It was
shown in \cite{BS} that this coincides with the usual Chern-Weil map $w_P$.

We show that this picture can be generalized to our more general context. While
it is possible to make sense of the above argument by replacing the manifold $M$
everywhere with a groupoid $\gm$, it is conceptually even simpler to broaden the
concept of ``maps'' between groupoids by allowing {\em generalized
homomorphisms}, which can be thought of as a groupoid version of smooth maps
between differential stacks. Under this framework, a principal $G$-bundle $P$
over a groupoid $\gm \toto \gm_0$ is then equivalent to a generalized
homomorphism from $\gm$ to $G$. Hence it induces a homomorphism on the level of
cohomology $H^*_{dR}(G_\com) \to H^*_{dR}(\gm_\com)$. Composing the Bott-Shulman
map with this map, one obtains a map $w: S(\Gg^*)^G \to H^*_{dR}(\gm_\com)$.
This is called the universal Chern-Weil map. We have

\vskip.2in

\noindent{\bf Theorem D.} {\it
The universal Chern-Weil map and the Chern-Weil map are equal.}

\vskip.2in

A connection on a principal $G$-bundle $P$ over a groupoid $\gm \toto \gm_0$ is
a pseudo-connection $\theta \in \Omega^1(P) \otimes \Gg$ such that
$\partial\theta \equal 0$. In contrast to pseudo-connections, the notion of
connections is well-defined for differential stacks in the sense that they may
pass to Morita equivalent groupoids in a natural fashion (Corollary
\ref{cor:3.15}). Unlike the case of manifolds, connections do not always exist
for every principal bundle over a groupoid $\gm \toto \gm_0$. However, we prove
that

\vskip.2in

\noindent{\bf Theorem E.} {\it
Any principal $G$-bundle over an orbifold always admits a connection.}

\vskip.2in

Whenever a connection exists, the Chern-Weil map admits a much simpler form. It
essentially results from applying a polynomial function to the curvature of the
connection exactly as in the manifold case. More precisely, consider the
subcomplex $(\Omega(\gm_0)^\gm, d)$ of the de Rham complex $(\Omega(\gm_0), d)$,
where $\Omega(\gm_0)^\gm \equal \{\omega \in \Omega(\gm_0)| \ \partial \omega
\equal 0\}$. Let $H^*_{dR}(\gm_0)^\gm$ be its cohomology. The inclusion $i:
\Omega(\gm_0)^\gm \to \Omega(\gm_{\com})$ is a chain map and therefore induces a
morphism $i: H^*_{dR}(\gm_0)^\gm \to H^*_{dR}(\gm_\com)$.

\vskip.2in

\noindent{\bf Theorem F.} {\it
Assume that the principal $G$-bundle $P \stackrel{\pi}{\to} \gm_0$ over the
groupoid $\gm \toto \gm_0$ admits a connection $\theta \in \Omega^1(P) \otimes
\Gg$. Then the following diagrams commute
    \[ \comtri{S(\Gg^*)^G}{z}{Z^*(\gm_0)^\gm}{z_\theta}{i}{Z^*(\gm_{\com})} \]
and
    \[ \comtri{S(\Gg^*)^G}{w}{H^*_{dR}(\gm_0)^{\gm}}{w_P}{i}
        {H^*_{dR}(\gm_{\com})} \]
where $z: S(\Gg^*)^G \to Z^*(\gm_0)^{\gm}$, and $w: S(\Gg^*)^G \to
H^*_{dR}(\gm_0)^{\gm}$ denote the usual Chern-Weil maps obtained by forgetting
the groupoid action. {\it I.e.}, $z(f) \equal f(\Omega)\in Z^*(P)^{basic} \cong
Z^*(\gm_0)$, $\forall f \in S(\Gg^*)^G$, where $\Omega \in \Omega^2(P) \otimes
\Gg$ is the curvature form.}

\vskip.2in

Therefore, as a consequence, when a flat connection exists, {\it i.e.} a
pseudo-connection whose total curvature vanishes, the Chern-Weil map vanishes
(except in degree $0$). Indeed, a flat connection resembles in many ways a flat
connection on usual principal bundles. In particular we prove

\vskip.2in

\noindent{\bf Theorem G.} {\it
A flat connection on a principal $G$-bundle $P \stackrel{\pi}{\to} \gm_0$ over a
groupoid $\gm \toto \gm_0$ induces a group homomorphism
    \[ \pi_1(\gm_\com,x) \to G, \]
where $\pi_1(\gm_\com,x)$ denotes the fundamental group of the groupoid $\gm$.}

\vskip.2in

Since our construction of the Chern-Weil map can be done purely algebraically,
following the idea of Cartan \cite{Cartan}, we carry out this construction for
any $G$-differential simplicial algebra, which we believe is of interest in its
own right. This is the content of Section 4. We present two equivariant
constructions. One is through fat realizations and the other is a functorial
approach in terms of simplicial algebras. The Chern-Weil map for principal
$G$-bundles over a groupoid is discussed in Section 5. Section 2 is preliminary
on principal $G$-bundles. We give several equivalent pictures of principal
$G$-bundles over a groupoid which should be of independent interest. Section 3
is devoted to the study of differential geometry of these  principal bundles
including connections, curvatures and holonomy maps.

Finally, a remark about notation. In this paper, $G$ always denotes a Lie group
and $\Gg$ its Lie algebra. By a  ``groupoid'', we always mean
a  ``Lie groupoid'' unless it is specified otherwise; the
source and target maps are denoted by $\so$ and $\ta$, respectively.

{\bf Acknowledgments.}
We would like to thank several institutions for their hospitality while work on
this project was being done: Penn State University (Tu), and RIMS and Erwin
Schr\"odinger Institute (Xu). We also  wish to thank Kai Behrend and Jim
Stasheff for useful discussions and comments.

\section{Principal $G$-bundles over groupoids}

\subsection{Simplicial $G$-bundles}

In this subsection, we recall some useful results concerning
simplicial $G$-bundles.

Recall  that a {\em simplicial set} $M_\com$ \cite{Dupont}
 is a sequence  of sets $(M_n)_{n \in \nn}$ together
with \degeneracy maps $ \epsilon_i^n: M_n \to M_{n-1}  ,i\equal 0,\dots,n$
and \face maps $ \eta_i^n: M_n \to M_{n+1}$, $i\equal 0,\dots,n$,  which
satisfy the simplicial identities:
$\epsilon_i^{n-1} \epsilon_j^{n}\equal  \epsilon_{j-1}^{n-1} \epsilon_i^{n}$
if $i<j$,
$\eta_i^{n+1} \eta_j^{n}\equal \eta_{j+1}^{n+1} \eta_i^{n} $ if $i \leq j$,
$\epsilon_i^{n+1} \eta_j^{n}\equal  \eta_{j-1}^{n-1} \epsilon_{i}^{n} $
if $i < j $, $\epsilon_i^{n+1} \eta_j^{n}\equal
\eta_{j}^{n-1} \epsilon_{i-1}^{n} $
if $i > j+1 $ and $\epsilon_j^{n+1} \eta_j^{n}\equal
\epsilon_{j+1}^{n+1} \eta_j^{n}\equal \id_{M_n} $.
A {\em simplicial manifold} is  a simplicial set  $M_\com\equal (M_n)_{n \in \nn}$,
where,  for every $n\in \nn$,  $M_n$ is a smooth manifold and all the \face 
and \degeneracy maps are smooth maps.

Associated to any Lie groupoid $\gm\toto \gm_0$,
there is a canonical simplicial manifold $\gm_\com$,
 constructed as follows \cite{Segal}.
 Set $ \Gamma_{n}\equal \{(\gamma_1,\gamma_2,\ldots,\gamma_n ) \in \Gamma^n 
| \so (\gamma_i)\equal \ta (\gamma_{i+1}),i\equal 1,\dots,n~-~1\}$,
 the manifold consisting of  all composable $n$-tuples, 
 and define the \degeneracy maps $\epsilon_{i}^n :\Gamma_{n}\to
\Gamma_{n-1}$ by, for $n > 1$
 \begin{eqnarray} 
\label{eq:simplicialite} 
 && \epsilon_{0}^n (\gamma_1  ,
 \gamma_2 , \ldots , \gamma_n)\equal  ( \gamma_2 , \ldots ,\gamma_n)
 \label{eq:face}\\
  && \epsilon_{n}^n (\gamma_1  , \gamma_2 ,
 \ldots , \gamma_n)\equal  (\gamma_1 , \ldots ,\gamma_{n-1})\\
 && \epsilon_{i}^n (\gamma_1  ,  \ldots,  \gamma_n)\equal 
 (\gamma_1 ,  \ldots, \gamma_i \gamma_{i+1} , \ldots ,\gamma_n), \ \  1\leq
i \leq n-1,
\end{eqnarray}
and for $n\equal 1$ by,
$ \epsilon_{0}^1 (\gamma)\equal \so (\gamma)$,
$\epsilon_{1}^1 (\gamma)\equal \ta (\gamma)$.
Also define the \face maps by $\eta_0^0\equal
\epsilon: \gm_0 \to \gm_1$ ($\epsilon$ being the unit map of the groupoid) and
$\eta^n_i :\Gamma_{n}\to \Gamma_{n+1}$ by:
\begin{eqnarray}
&& \eta^n_0 (\gamma_1 ,\ldots , \gamma_n)\equal 
(\ta ( \gamma_1), \gamma_1 , \ldots , \gamma_n) \label{eq:deg}\\
&& \eta^n_i (\gamma_1 ,\ldots , \gamma_n)\equal 
(\gamma_1 , \ldots ,\gamma_{i}, \so (\gamma_{i})
,\gamma_{i+1}, \ldots , \gamma_n), \ \ 1\leq i\leq n.
\end{eqnarray}

We now recall the notion of  simplicial $G$-bundles
over a simplicial manifold  \cite{Dupont}.

\begin{defn}
\label{def:simplicial}
A principal $G$-bundle over a simplicial manifold
$M_\com: \equal (M_n)_{n \in \nn}$ is a simplicial manifold $P_\com:
\equal (P_n)_{n \in \nn}$ such that
\begin{itemize}
\item     for every  $n \in \nn$, $ P_n$ is a principal $G$-bundle over $M_n$, 
and
\item    the \face and \degeneracy maps are morphisms
of principal $G$-bundles. 
\end{itemize}
\end{defn}

\subsection{Principal $G$-bundles over a groupoid} \label{subsection:principal}

We recall in this subsection
 the definition of a principal $G$-bundle
over a groupoid  \cite{TXL,Moe},
 and give some basic examples.

\begin{defn} \label{def:principal} 
 A principal $G$-bundle over a  groupoid
 $\Gamma\rightrightarrows \Gamma_0$ consists of  a principal right
$G$-bundle $P \stackrel{\pi}{\to} \Gamma_0$ over the manifold $\Gamma_0$
 such that

(i) there is a map $\sigma : Q\to P$, where $Q$ is the fibered
product $Q\equal \gm \times_{\so, \Gamma_0, \pi} P$.
We write $\sigma(\gamma, p)\equal \gamma \cdot p$.
 This map is subject to the constraints

(ii) for all $p\in P$ we have  $\pi(p) \cdot p\equal p$;

(iii)  for all $p\in P$ and all $\gamma_1,\gamma_2\in\Gamma$ such that
$\pi(p)\equal \so(\gamma_2)$ and $\ta(\gamma_2)\equal \so (\gamma_1)$ we have
$  \gamma_1 \cdot (\gamma_2  \cdot p)\equal  (\gamma_1  \gamma_2)  \cdot p$;

(iv) for all $p\in P$, 
$g\in G$ and $\gamma \in \gm$, such that $\so(\gamma ) \equal \pi(p)$, we have
$ (\gamma \cdot p) \cdot g \equal \gamma \cdot (p \cdot g )$;
\end{defn}

Axioms (i)-(iii) simply mean that the groupoid $\Gamma \toto \Gamma_0$
acts on $P  \stackrel{\pi}{\to} \Gamma_0 $.
Axiom (iv) means that this action commutes with the $G$-action.



%


\begin{numex}
\label{ex:equivariant}
Let $\gm$ be a transformation groupoid $H\times M\toto M$,
where $H$ acts on $M$ from the left.
Then a principal $G$-bundle over $\gm$ corresponds exactly
to an $H$-equivariant principal (right) $G$-bundle over $M$.
\end{numex}

The following proposition 
 clarifies the relation between principal $G$-bundles
over a groupoid $\Gamma \toto \Gamma_0 $ and 
principal $G$-bundles over its  corresponding
simplicial manifold $\Gamma_\com$.

\begin{prop} \label{prop:11corresp}
Let $\Gamma \toto \Gamma_0 $ be a Lie groupoid.
There is an  equivalence   of categories 
between the category of $G$-bundles
over the groupoid $\Gamma \toto \Gamma_0$ and the category of simplicial
$G$-bundles over the simplicial manifold $ \Gamma_\com$.
 \end{prop}

This proposition is an immediate consequence of
Lemma \ref{lem:G-morphism} and Lemma  \ref{lem:eq2} below.

First, let us introduce some terminology.
For a given Lie groupoid $\Gamma \toto \Gamma_0 $,
by a {\em principal $G$-groupoid over } $\Gamma \toto \Gamma_0 $,
we mean   a Lie groupoid $Q\toto P$  together with a 
 groupoid morphism from $Q$ to $\Gamma$:
 \begin{equation}
\label{eq:Gpoid}
  \begin{array}{ccccc}
 G & \to &  Q &\stackrel{\pi}{\to} & \Gamma     \\
   & &  \downdownarrows       &    &    \downdownarrows         \\
 G &\to  &  P &\stackrel{\pi}{\to} & \Gamma_0   \\
\end{array}
\end{equation}
  such that  both  $G\to Q\stackrel{\pi}{\to}  \Gamma $ 
 and  $G\to P\stackrel{\pi}{\to}  \Gamma_0 $
are principal $G$-bundles and the \source and \target maps
are morphisms of principal $G$-bundles. There is an obvious notion
of morphisms for principal $G$-groupoids over a 
given groupoid $\gm\toto \gm_0$ so that one obtains a 
category.                               

\begin{lem}
\label{lem:G-morphism}
The category  of  principal $G$-bundles over
the groupoid $\Gamma \toto \Gamma_0$ is equivalent   to the
category   of   principal $G$-groupoids over $\Gamma \toto \Gamma_0$.
\end{lem}     
\begin{pf}
 Associated to any principal $G$-bundle  $ P\stackrel{\pi}{\to} \Gamma_0$
 over $\gm \toto \gm_0$, there is a groupoid $Q\toto P$,   namely,
the {\em transformation groupoid},
which is defined as follows.
Let   $Q\equal\gm\times_{\so, \Gamma_0,\pi}P$,
the \source and \target maps are, respectively,
$\so (\gamma , p)\equal  p$, $\ta (\gamma , p)\equal \gamma \cdot  p$,
  and the multiplication is\begin{equation}
\label{eq:transformation}
(\gamma_1 ,q )\cdot (\gamma_2 , p)\equal (\gamma_1 \gamma_2 , p),
\ \ \ \ \mbox{where }
q\equal \gamma_2 \cdot p.
\end{equation}
Let $\pi : Q\to \gm$ be the projection map. It is simple to see
that $\pi$ is a groupoid homomorphism and one obtains  a commutative diagram
(\ref{eq:Gpoid}). Hence $Q$ is a principal $G$-groupoid over
$\gm$. 
   

Conversely, assume that $Q$ is a principal $G$-groupoid over
$\gm$. 
Define a  $\gm$-action  on $ P  \stackrel{\pi}{\to}  \Gamma_0$
 by $ \gamma \cdot p\equal  \ta(q)$,  where $q \in Q$
is the unique element on the fiber  $\pi^{-1} ( \gamma )$ 
satisfying  $\so(q)\equal p$.
This  endows $ P \stackrel{\pi}{\to} \Gamma_0$ with a structure 
of principal $G$-bundle over $\gm \toto \gm_0 $ such that $ Q \toto P$  is 
its corresponding transformation groupoid.

We leave the reader to  check the functoriality of the above construction.
\end{pf}

\begin{lem}
\label{lem:groupoidinduit}
Let  $\Gamma \toto \gm_0$ be a Lie  groupoid.
Assume that   $Q$ and $P$ are   principal $G$-bundles  over   $\Gamma$ and
$\gm_0$,  respectively,
 which admit   two $G$-bundle maps $\so,  \ta : Q\to P$ such
that  diagram (\ref{eq:Gpoid}) commutes.
 Then $Q$ admits  a structure of principal
 $G$-groupoid over $ \gm \toto \gm_0$.
\end{lem}
\begin{pf} 
We take $\so:Q \to  P$ and $ \ta:Q \to P$
as the  \source and \target maps.
The product $q_1\cdot q_2 $ of two elements $q_1$ and $q_2$
 with $\so(q_1)\equal \ta(q_2)$
is  defined to be 
the unique element on  the $\pi$-fiber over $ \pi(q_1) \pi(q_2) \in \Gamma$
whose  \source equals to  $\so(q_2)$.  For any $p\in P$, define the unit
$\epsilon (p)\in Q$ to be the unique element in 
$\pi^{-1}(\epsilon (\pi (p))$ such that $\so (\epsilon (p) )\equal p$.
For any $q \in Q$,
the inverse $q^{-1}$ is defined  to be the unique element
on the  $\pi$-fiber over $\pi(q)^{-1} \in \Gamma$
such that $\so (q^{-1})\equal \ta (q)$.

Since for any composable  triple $(q_1,q_2,q_3)$,
 both $q_1\cdot (q_2 \cdot q_3) $
and $(q_1\cdot q_2)\cdot q_3 $ 
are mapped  to  $\pi(q_1) \pi(q_2)\pi(q_3) \in \Gamma $ under $\pi$, and
to $ \so(q_3)$ under $\so$, hence $q_1\cdot (q_2 \cdot q_3) \equal (q_1\cdot q_2)\cdot q_3$.
Therefore the multiplication
 is associative. The rest of the axioms can be verified 
in a similar fashion.
\end{pf}

\begin{lem}\label{lem:eq2}
The category of    principal $G$-groupoids over $\Gamma \toto \Gamma_0$
is equivalent  to the category   
of simplicial
$G$-bundles over the simplicial manifold $ \Gamma_\com$.    
\end{lem} 
\begin{pf}
 Given a principal $G$-groupoid $Q\toto P$
 over $\Gamma \toto \Gamma_0 $,
$ Q_\com$ is clearly 
a principal $G$-bundle over the simplicial manifold $\Gamma_\com $.

Conversely, assume that $ R_{\com} :\equal (R_n )_{n\in \nn }$ is
 a simplicial  $G$-bundle over $\Gamma_{\com}$.
Let $Q\equal R_1 $ and $P\equal R_0 $.
By Definition \ref{def:simplicial},
 we have the    commutative diagram  as in Eq. (\ref{eq:Gpoid}).
Therefore, according to Lemma \ref{lem:groupoidinduit}, $Q\toto P$
is  endowed with  a groupoid structure, which  is a principal $G$-groupoid
over $\gm \toto \gm_0$.


To complete the proof, it suffices  to  show that 
the associated simplicial manifold $ Q_\com$ is
isomorphic to $R_\com $ as simplicial $G$-bundles over $ \Gamma_\com$.
By abuse of notation, we use $\epsilon_i^n$ to denote 
the \degeneracy maps of both  $ R_\com$ and $ Q_{\com}$. 
We denote by $ p_n$   the maps from $ Q_n $ or   $R_n$ to $P$
defined by $ p_n\equal  \epsilon_0^1 \smalcirc \dots \smalcirc \epsilon_0^n $.
Define $\phi_n: R_n \to Q_n$  by
$\phi_n (r) \equal q $, 
 where  $q \in Q_n$ is the unique element   in the
fiber  of $\pi :Q_n \to \gm_n$
 over  $ \pi(r ) \in \gm_n $ with $ p_n(q)\equal p_n(r)$. It is simple
to see that $(\phi_n)_{n \in \nn}$ 
is  a simplicial map and is therefore an isomorphism of 
principal $G$-bundle over $\gm_\com$. 
 \end{pf}

\subsection{Generalized homomorphisms}

In this subsection, we  will give another definition
 of  $G$-bundles over a groupoid
using the notion of generalized
 homomorphisms \cite{hae84,hilsum-skandalis87,TXL}.
Let us  recall  its definition below.

\begin{defn}\label{def:generalized morphism}
A generalized groupoid homomorphism $\phi: \equal (Z,\sigma,\tau)$
from $\Gamma \toto \Gamma_0$
to $H \toto H_0$ is given by a manifold $Z$, 
two smooth maps   $\Gamma_0 \stackrel{\tau }{\leftarrow}
Z \stackrel{\sigma}{\rightarrow} H_0$,
a left action of $\Gamma$ with respect to $\tau$,
a right action of $H$ with respect to $\sigma$, such that the two
actions commute, and  $Z$ is an  $H$-principal
bundle over $\Gamma_0$.
\end{defn}

In the sequel, we will use both notations $(Z,\sigma,\tau)$
and $\Gamma_0 \stackrel{\tau }{\leftarrow}
Z \stackrel{\sigma}{\rightarrow} H_0$ interchangeably to denote
a generalized  homomorphism.

Recall that given a Lie groupoid $H\toto H_0$, by
an (right)  {\em $H$-principal bundle over $M$},
 we mean an  (right) $H$-space $P\to H_0$ and a surjective
submersition $\pi: P\to M$ such that
 for all $p,p'\in P$, such that $\pi(p)=\pi(p')$, there
exists a unique $\gamma\in  H$, such that $p\cdot\gamma$ is defined
and $p\cdot\gamma=p'$.

The following proposition lists  several useful equivalent definitions 
of  (right) $H$-principal
bundles, or (right) $H$-torsors over a manifold $M$.

\begin{prop}
(\cite[\S 2]{hae84}, \cite[Def. 1.1]{hilsum-skandalis87})
Let  $H \toto H_0$ be a Lie groupoid.
The following statements are equivalent.

\begin{enumerate}
\item $Z\stackrel{\sigma}{\rightarrow} H_0$ is
an  (right) $H$-principal bundle over $M$;
\item $Z\stackrel{\sigma}{\rightarrow} H_0$ is an right
$H$-space, where the $H$-action is free and proper so
that $M\cong Z/H$ (note that the projection map
$Z\to M$ is always a surjective submersion).
\item  $Z\stackrel{\sigma}{\rightarrow} H_0$ is  an right
$H$-space  with a smooth map $\pi: Z\to M$ 
such that at every point $x\in M$ there is
a local section  $s: U\to Z$ of $\pi$
satisfying $\pi^{-1}(U) \cong U\times_{\sigma\smalcirc s, H_0, \ta} H$.
\end{enumerate}
\end{prop}

If $f\colon \Gamma\to H$ is a 
strict homomorphism
({\em i.e.} a smooth map satisfying the property  that 
$f(\gamma_1 \gamma_2)\equal f(\gamma_1)f(\gamma_2)$),  then
$Z_f\equal \Gamma_0 \times_{f,H_0,\ta} H $, with $\tau(x,h)\equal x $,
$\sigma(x,h) \equal \so(h)$, and the actions
$\gamma\cdot(x,h)\equal \big( \ta (\gamma),f(\gamma)h\big)$ and
$(x,h)\cdot h'\equal (x,h h')$, defines a generalized homomorphism
from $\Gamma \toto \Gamma_0$ to $H \toto H_0$.

Two generalized homomorphisms
$\phi_1:\equal (Z,\sigma,\tau)$ and $\phi_2:\equal (Z',\sigma',\tau')$
from $ \Gamma$ to $H$ are said to be {\em equivalent} if 
 there is a $\gm$-$H$-equivariant diffeomorphism $\phi: Z\to Z'$.

Generalized homomorphisms can be composed as follows.
If  $\phi: \Gamma_0 \stackrel{\tau }{\leftarrow} Z 
\stackrel{\sigma}{\rightarrow} H_0$ and
$\phi' : H_0 \stackrel{\tau' }{\leftarrow} Z' 
\stackrel{\sigma'}{\rightarrow} R_0$ 
are generalized homomorphisms from
$\gm\toto \gm_0 $ to $H\toto H_0$, and from
$H \toto H_0$ to $R\toto R_0$,
 respectively, then the composition $ \phi \smalcirc \phi':
 \gm_0 \leftarrow Z'' \rightarrow R_0$ defined by 
$$Z'' \equal Z\times_H Z':\equal {(Z\times_{\sigma,H_0,\tau'}Z')}
/_{(z,z')\sim(zh,h^{-1}z')}$$
is a generalized homomorphism
from $\Gamma\toto \Gamma_0  $ to $R\toto R_0 $.
 Moreover, the composition
of generalized homomorphisms is associative. Thus one obtains  a
category ${\mathcal{G}}$ whose objects are Lie groupoids and morphisms
are generalized homomorphisms \cite{hilsum-skandalis87,TXL}.
 There is a functor ${\mathcal{G}}_s\to {\mathcal{G}}$,
where ${\mathcal{G}}_s$ is the category of Lie groupoids  with
strict homomorphisms, given by $f\mapsto Z_f$ as described above.
The  isomorphisms in the category ${\mathcal{G}}$ are called
 {\em Morita equivalences} \cite{muhly-renault-williams87, Xu:90}.
Indeed it is simple to see that any generalized homomorphism
can be decomposed as the composition of a Morita equivalence
with a  strict homomorphism  \cite{TXL}. Let us
briefly  recall  this construction below.

Let $\phi:\equal (Z,\sigma,\tau)$ be a generalized homomorphism
from $\Gamma' \toto \Gamma_0' $ to $\Gamma \toto \Gamma_0 $.
 We denote  by $\Gamma'[Z]\toto Z$ the pull-back of
 $\Gamma'\toto \gm_0'$ via the
surjective submersion $Z\stackrel{\tau }{\to}\Gamma_0'$, {\em {\it i.e.}}
 the groupoid $Z\times_{\tau,\ta}\Gamma'\times_{\so,\tau} Z$ 
with the  multiplication law
$(z_1,\gamma_1',z_2)(z_2,\gamma_2',z_3)\equal (z_1,\gamma_1'\gamma_2',z_3)$.
Then  the projection map from
 $\Gamma'[Z]\toto Z $ to $\Gamma'\toto \gm_0'$, denoted by $\tau$,
by abuse of notation,  
\begin{equation}\label{eq:deftau}  \tau(z_1,\gamma',z_2)
\equal  \gamma' \end{equation}
is a strict homomorphism, which is indeed a Morita equivalence \cite{BX}.
On the other hand  the  map from
 $\Gamma'[Z]\toto Z $ to $\Gamma\toto \gm_0$,
\begin{equation}\label{eq:defsigma} 
 \sigma(z_1,\gamma',z_2)\equal  \gamma, 
\end{equation}
where  $\gamma  \in \gm$ is the  unique element such that 
$\gamma' \cdot z_2\equal z_1\cdot \gamma$, is  a strict homomorphism.
Thus we have proved the following

\begin{lem}\label{lem:gener morphism strict}
Let $\phi\equal (Z,\sigma,\tau)$ be  a  generalized homomorphism 
from  $\Gamma' \toto \Gamma_0' $ to $\Gamma \toto \Gamma_0 $.
Then $\phi$  is   the composition of the Morita equivalence
$\Gamma' \sim_{Morita} \gm' [Z]$ with the strict homomorphism $\sigma: 
\gm' [Z]\to \gm$.
\end{lem}

The following proposition follows  immediately from the definition
of  generalized homomorphisms.

\begin{prop} 
\label{prop:11corr}
 The category of $G$-principal bundles over a groupoid
 $ \Gamma \toto \Gamma_0$ is equivalent to the category of
generalized homomorphisms from  $ \Gamma \toto \Gamma_0$ to $G\toto \cdot$.
\end{prop}

As a consequence, we have the following

\begin{cor}
 If  $\gm'$ and  $\gm$ are Morita equivalent groupoids,
there is an equivalence 
of categories of principal $G$-bundles  over $\gm'$ and  $\gm$.               
\end{cor}

As another consequence, principal $G$-bundles over groupoids 
can be pulled back via generalized homomorphisms.
We describe this pull-back construction explicitly in the following: 

\begin{prop}
\label{prop:pullback}
Let $P$ be a principal $G$-bundle over the groupoid 
$\gm$ and $\phi :  \Gamma'_0 \stackrel{\tau }{\leftarrow} Z
\stackrel{\sigma}{\rightarrow} \gm_0$ a generalized
homomorphism from $\gm'$ to $\gm$. Then
 the pull-back  $G$-bundle $ P' \stackrel{\pi'}{\to} \gm_0'$
can be described as follows. 
\begin{itemize}
\item $P'\equal (Z\times_{\Gamma_0 }P)/\gm $.  We  denote by
$\overline{(z,p)} \in P'$  the  class in $P'$
corresponding to $ (z,p) \in Z \times_{\Gamma_0 }P $; 
\item   $G$  acts on $P'$ by
  $$ \overline{(z,p)}\cdot g \equal
  \overline{(z,p \cdot g)}, \hspace{1cm}  \forall g \in G;$$
\item $\pi' : P'\to \gm_0' $ is given by
 $$ \pi'\big( \overline{(z,p)}\big) \equal \tau(z); $$
\item the  $\gm'$-action on $P'\to \gm_0' $ is given by 
  $$ \gamma'  \cdot \overline{(z,p)} \equal  \overline{(\gamma' \cdot z,p)},
 \hspace{1cm}  \forall \gamma' \in \gm'  \ \ \ 
 \mbox{ with $\so(\gamma')\equal \tau(z)$}.$$
\end{itemize}
\end{prop}

Sometimes it is  useful to   consider
a  generalized homomorphism  $\phi:\equal (P,\sigma,\pi) $
defined by a principal $G$-bundle $P \stackrel{\pi}{\to} \gm_0$ 
over a groupoid $\gm \toto \gm_0$
as the   composition of a  Morita equivalence with a  strict homomorphism.
Below we present two equivalent pictures of such compositions.

First  consider the semi-direct  product groupoid 
${Q}': \ Q\times G\toto P$ (using the right $G$-action
 on  the groupoid $Q\toto P$ as introduced in the proof of 
Lemma \ref{lem:G-morphism}).
Here the \source and \target maps,
and the multiplication are
\begin{equation}\label{eq:productQ'}
\so(q, g)\equal \so (q)g,     \, \, \, \ta(q, g)\equal \ta (q), \, \, \,
(q_1 , g_1 )\cdot (q_2 , g_2 )\equal (q_1 \cdot (q_2 g_1^{-1})\ ,\ 
 g_1 g_2).
\end{equation}
It is simple to check that this groupoid  is Morita equivalent to
$\gm\toto \gm_0$. In fact
the map $(\gamma,p,g)\mapsto (\gamma p,\gamma,p g)$
is a groupoid  isomorphism from $Q'\toto P$ onto the pull-back groupoid
of $\gm \toto \gm_0$ via the surjective submersion
$P\stackrel{\pi}{\to} \gm_0$.
Moreover the natural  projection on the second factor
\begin{equation}
\label{eq:Q}
\begin{array}{ccc}
 Q'\equal Q\times G&\stackrel{pr}{\to} &G \\
  \downdownarrows & &  \downdownarrows\\
 P &\to &\cdot \\
\end{array}
\end{equation}
defines a strict  groupoid homomorphism.

Alternatively, consider 
the gauge groupoid  $\frac{P \times P}{G} \toto \Gamma_0$.
We denote  by $\overline{(p_1, p_2 )}$ the class 
 in $\frac{P \times P}{G}$ corresponding to $(p_1, p_2)\in P\times P$.
A strict homomorphism $ \rho :\Gamma \to
 \frac{P \times P}{G}$ is 
defined by $\gamma \mapsto
\overline{(\gamma p, p)}$,  where $p\in P$ is any element satisfying
$\pi(p)\equal \so (\gamma)$.
 Thus we obtain the following groupoid homomorphism:
\begin{equation}
\label{eq:Q1}
\xymatrix{ 
\Gamma \ar@<-.5ex>[d]\ar@<.5ex>[d]\ar[r]^\rho &
\frac{P \times P}{G} \ar@<-.5ex>[d]\ar@<.5ex>[d]\\
\gm_0 \ar[r] & \gm_0}
\end{equation}
                                         
Since any transitive groupoid is Morita equivalent  to its isotropy group,
$\frac{P\times P}{G}\toto M $ is Morita equivalent to $G \toto \cdot$.

\begin{prop}
Let $P \stackrel{\pi}{\to} \gm_0$ be a principal $G$-bundle
over a groupoid $\gm \toto \gm_0$. Then the generalized homomorphism 
 $\Gamma_0 \stackrel{\pi }{\leftarrow} P \rightarrow
 \cdot $ from $ \Gamma \toto \Gamma_0$ to $G\toto \cdot$ is  equal to
 the composition of $\gm\sim_{\mathrm{Morita}} Q'$ and $\pr: Q'\to G$,
as well as the composition of $\rho: \gm \to     \frac{P \times P}{G}$ 
and $\frac{P \times P}{G}\sim_{\mathrm{Morita}} G$.
\end{prop} 
\begin{pf}
This follows from  a direct verification,
 and is left to the reader.
\end{pf}

We end this section by  summarizing
 Proposition \ref{prop:11corresp}
 and  Proposition \ref{prop:11corr} in the following:

\begin{them} 
Let $\Gamma \toto \Gamma_0$ be  a Lie groupoid.
The following categories are equivalent:
\begin{itemize}
\item the category of  principal $G$-bundles over  $\Gamma \toto \Gamma_0$;
\item the category of  principal $G$-bundles
over the simplicial manifold $\Gamma_{\com} $; and
\item the category of   generalized homomorphisms from
 $\Gamma \toto \Gamma_0$ to $G \toto \cdot$.
\end{itemize}
\end{them}

\section{Pseudo-connections and connections}

The purpose of this section is to study differential
geometry of a principal $G$-bundle over a Lie  groupoid
 and in particular connections and curvatures.

\subsection{De Rham cohomology}
\label{section:deRham}
Considering a Lie groupoid $\gm\toto\gm_0$ as a generalization
of a manifold (indeed it defines a differential stack \cite{BX1}), 
one can also speak about de Rham cohomology, which is
given by the following double complex $\Omega^\com(\Gamma_\com)$:

\begin{equation}
\label{eq:DeRham}
\xymatrix{
\cdots&\cdots&\cdots&\\
\Omega^1(\gm_0)\ar[u]^d\ar[r]^\partial &\Omega^1(\gm_1)\ar[u]^d\ar[r]^\partial
&\Omega^1(\gm_2)\ar[u]^d\ar[r]^\partial&\cdots\\
\Omega^0(\gm_0)\ar[u]^d\ar[r]^\partial&\Omega^0(\gm_1)\ar[u]^d\ar[r]^\partial
&\Omega^0(\gm_2)\ar[u]^d\ar[r]^\partial&\cdots
}
\end{equation}
Its boundary maps are $d:
\Omega^{k}( \gm_p ) \to \Omega^{k+1}( \gm_p )$, the usual exterior
differential,  and $\partial
:\Omega^{k}( \gm_{p} ) \to \Omega^{k}( \gm_{p+1} )$,  the alternating
sum of the pull-back  of the face maps  of the corresponding
simplicial manifolds $\gm_\com$.
We denote the total differential by 
$$\delta\equal (-1)^pd+\del .$$
The cohomology groups of the total complex 
$H_{dR}^k(\Gamma_\com)\equal H^k\big(\Omega^\com(\Gamma_\com)\big)$
are called the \emph{de~Rham cohomology} 
groups of $\Gamma\toto \gm_0$.              
Note that $H_{dR}^* (\Gamma_\com ) $ is a super-commutative algebra
with respect to the cup-product \cite{Dupont}
 defined  as follows. For any
 $\omega_1 \in \Omega^k(\Gamma_{m})$ and $\omega_2 \in \Omega^l(\Gamma_{n})$,
   let $\omega_1 \vee \omega_2 \in \Omega^{k+l} (\gm_{m+n} )$ be the
differential form given by 
\begin{equation}
\label{eq:cup}
(\omega_1 \vee \omega_2 ) (\gamma_1,\dots,\gamma_{m+n})=
(-1)^{kn} p_1^*\omega_1 \wedge p_2^* \omega_2 ,
\end{equation}
where $p_1:\Gamma_{m+n}\to \Gamma_m$ is defined by
$$ \left\{ \begin{array}{cccl} p_1(\gamma_1,\dots,\gamma_{m+n}) & 
\equal&(\gamma_1,\dots,\gamma_m) & \mbox{if $m \geq 1$} \\
 p_1(\gamma_1,\dots,\gamma_{n})&\equal&\alpha(\gamma_1)& \mbox{if $m \equal0
 $ and $n \geq 1$} \\  p_1  &\equal &Id  & \mbox{if $m\equal n 
 \equal 0$}   \end{array} \right. $$
and  $p_2:\Gamma_{m+n}\to \Gamma_n$ is the map
$$ \left\{ \begin{array}{cccl} p_2(\gamma_{1},\dots,\gamma_{m+n}) &  \equal 
&(\gamma_{m+1},\dots,\gamma_{m+n} ) & \mbox{if $n \geq 1$} \\
 p_2(\gamma_1,\dots,\gamma_{m})& \equal &\beta(\gamma_m)& \mbox{if $n  \equal 0 $}  \\ 
 p_2  & \equal  &Id  & \mbox{if $m \equal n \equal 0$}
  \end{array} \right. $$
Eq. \eqref{eq:cup} can be rewritten in
terms of face maps into a  more compact form:

\begin{equation}\label{eq:cup-productdR} 
 \omega_1 \vee  \omega_2  \equal  (-1)^{kn}
(\epsilon_{m+n}^{m+n})^*  \cdots  (\epsilon_{m+i}^{m+i})^* 
  \cdots   
(\epsilon_{m+1}^{m+1})^*(\omega_1) \wedge (\epsilon_{0}^{n+m})^*
  \cdots   (\epsilon_{0}^{n+i})^*  \cdots  
  (\epsilon_{0}^{n+1})^*(\omega_2).
\end{equation}

One can prove that $\vee$ induces a graded ring structure
on the de Rham cohomology groups $ H^*_{dR}(\gm_\com)$.
An important property of the de Rham cohomology is that it is
 functorial with respect to generalized homomorphisms. More precisely, 
a generalized homomorphism $\phi$ from  $ \Gamma' \toto \Gamma_0'$
to $ \Gamma \toto \Gamma_0$
induces a homomorphism 
 $\phi^*:H^*_{dR}(\gm_\com) \to H^*_{dR}(\Gamma'_\com)$.
We  recall its  construction below.

First, note that  a strict homomorphism $\psi$ of Lie groupoids
induces a  natural chain map between  their de Rham
double complexes, and therefore a morphism $\psi^*$ of de Rham cohomology groups.
 Now given a generalized  homomorphism  $\phi:\equal (Z,\sigma,\tau)$
 from $ \Gamma' \toto \Gamma'_0$ to  $ \Gamma \toto \Gamma_0$,
 let  $ \Gamma'[Z] \toto Z $ be the pull-back
 groupoid of $ \Gamma' \toto \Gamma_0' $ under  the
surjective  submersion $ \tau: Z \to \Gamma'_0$.
Denote again by $\tau$  the strict homomorphism 
from $ \Gamma'[Z] \toto Z $
to  $ \Gamma' \toto \Gamma_0' $,  and 
the strict homomorphism  from $ \Gamma'[Z] \toto Z $ to 
$ \Gamma \toto \Gamma_0 $ as in Eq.~(\ref{eq:deftau})
and (\ref{eq:defsigma}),  respectively.  Then

\begin{lem}
\cite{BX, tu05}
\label{lem:moritamap} The homomorphism 
 $ \tau^*: H^*_{dR}(\Gamma'_\com) \iso H^*_{dR}(\Gamma' [Z]_\com)  $
is an isomorphism.
\end{lem}

Define $\phi^*: \ H^*_{dR}(\gm_\com) \to H^*_{dR}(\Gamma_\com' )$ 
by

\begin{equation} \label{eq:defpullback}
\phi^* \equal  (\tau^*)^{-1} \smalcirc \sigma^* .
 \end{equation}
It is simple to check that two equivalent generalized
homomorphisms induce the same  morphism
for the de Rham cohomology groups. Hence $\phi$ is well-defined.
 In particular, if $ \Gamma' \toto \Gamma'_0$ and   $ \Gamma \toto \Gamma_0$
 are Morita equivalent groupoids, then $ \ H^*_{dR}(\gm_\com)  \iso
  H^*_{dR}(\Gamma'_\com )$.

\subsection{Pseudo-connections and pseudo-curvatures}

We are now ready to introduce the following

\begin{defn}
Let $P\to \gm_0$ be a  principal $G$-bundle over a Lie groupoid
$\gm\toto\gm_0$.
\begin{itemize}
\item A  pseudo-connection is a connection
 one-form $\theta  \in \Omega^1 (P)\otimes \Gg$ of the
$G$-bundle $P\to \gm_0$ (ignoring the groupoid action);
\item the total  pseudo-curvature $\Omega_{total}\in \Omega^2( Q_\com )
\otimes \Gg$ is defined by
\begin{equation}\label{eq:defomegatotal}
\Omega_{total} \equal  \delta \theta +\frac{1}{2}[\theta,\theta]\equal 
 \partial \theta +\Omega,
\end{equation}
where $\Omega \equal d\theta + \frac{1}{2}[\theta, \theta] \in
\Omega^2 (P)\otimes \Gg$ is the curvature form corresponding to $\theta$.
\end{itemize}
\end{defn}
The total pseudo-curvature $ \Omega_{total} $ consists of 
 two terms: the term $\partial \theta :\equal \so^{*}\theta-\ta^{*}\theta $
  is a  $\Gg$-valued $1$-form on $Q$,  and
 the term  $\Omega$  is a   $\Gg$-valued $2$-form on $ P$.
Both terms have a total degree  $2$
in  the double complex $\Omega^\com(Q_{\com}) \otimes \Gg$.

Introduce a  bracket of degree 
$-1$ on $\Omega^\com (Q_{\com}) \otimes \Gg$ by
   \begin{equation} \label{eq:defbra} [\omega_1 \otimes X_1,\omega_2 \otimes X_2]\equal 
(\omega_1 \vee \omega_2 ) \otimes [X_1,X_2] ,
\end{equation} 
where $\omega_1,\omega_2 \in \Omega^\com (Q_\com)$,
and  $X_1,X_2 \in \Gg$.
Note that this bracket does not satisfy the graded-Jacobi 
identity,  and is neither graded skew-symmetric.

\begin{prop}
The total pseudo-curvature satisfies the Bianchi identity:
\begin{equation}
\label{eq:11}
   \delta \Omega_{total} \equal \half
 ( [\Omega_{total}, \theta]-[\theta, \Omega_{total}] ). 
 \end{equation}
\end{prop}
\begin{pf}
 Applying the total differential $\delta$ to Eq.~(\ref{eq:defomegatotal}),
we obtain
   \be 
  \delta \Omega_{total}
&\equal & \frac{1}{2}\delta  ([\theta, \theta])    \\ 
 && \equal  \frac{1}{2} ([\delta \theta, \theta ]
-  [ \theta,\delta \theta ]) \\
 && \equal  \frac{1}{2} ([\Omega_{total}, \theta ] -
[ \theta,\Omega_{total} ]) 
+\frac{1}{4} ([[\theta,\theta],\theta]- [\theta, [\theta,\theta]]).
  \ee
Note that $\theta $ is an element in $ \Omega^1  (P) \otimes \Gg$,
and  the restriction of the cup-product to $\Omega (P)$ is simply
 the wedge product.
It follows from the Jacobi identity that
$[[\theta,\theta],\theta] \equal [\theta,[\theta,\theta]]\equal 0$  \cite{KOb}.
Therefore, Eq. (\ref{eq:11})
follows. 
\end{pf}

\begin{numrk}
We end this subsection by presenting a geometric
 interpretation  of $\partial \theta$.

Note that any  smooth path $p(t) $ in $P$ induces a smooth path
$g_p (t)$ in $G$ as follows.
By  $\tilde{p}(t) \in P   $, we denote the   horizontal
lift of the   projected path $\pi (p(t) )$ starting at
$p(0)$. Then  the path $ g_p (t) \in G $ is defined by the relation
$p(t)\equal  \tilde{p}(t) \cdot g_{p}(t) $. 
It is clear that $ g_p (t) $ satisfies the relation
   $\dot{g}_p (0) \equal  \dot{p}(0) \per \theta $.
Let $\delta_q \in T_qQ$ be any tangent  vector
and $q(t)$
 a path in $Q$ with $\dot{q}(0)\equal \delta_q$.
Set $a(t)\equal  \ta (q(t)) $ and $b(t)\equal  \so (q(t)) $.
Define  $g (t)\in G$ by 
$ g(t) \equal  g_a^{-1}(t) g_b(t)$. Then we have
\begin{equation}
\label{eq:interp} 
g(0)\equal 1, \ \ {\mbox{and}} \ \ \dot{g}(0)\equal \delta_q  \per \partial \theta.
\end{equation}
\end{numrk}


\subsection{Connections}

\begin{defn}
\label{def:connection}
Let $P\to \gm_0$ be a  principal $G$-bundle over a Lie groupoid
$\gm\toto\gm_0$.
\begin{itemize}
\item A   pseudo-connection
  $\theta  \in \Omega^1 (P)\otimes \Gg$ is called a connection if
 $\partial \theta \equal 0$;
\item if $\theta  \in \Omega^1 (P)\otimes \Gg$ is a connection, then  
$\Omega_{total}\equal \Omega\in \Omega^2 (P)\otimes \Gg$ is called
the curvature; 
\item a connection is said to be flat if its curvature vanishes.
\end{itemize}   
\end{defn}

In what follows, we investigate the criteria for a  connection to exist.
As we see below, this
 imposes a strong assumption on the  $G$-bundles.

By $U_{loc}(\Gamma)$ we denote  the pseudo-group of local bisections  of the
groupoid $\gm\toto \gm_0$. It is simple to see that
there is a local action of
 $U_{loc}(\Gamma) $  on  $P$ 
preserving the  $G$-bundle structure. In other words,
there is a group homomorphism
$$\phi : U_{loc}(\Gamma) \lon \Diff_{loc}(P), $$
where $\Diff_{loc}(P)$ denotes the pseudo-group
of local  automorphisms of the principal  $G$-bundle $P$. The map  $\phi$
 can be defined as follows. 
For any $\LL  \in U_{loc}(\Gamma)$, $\phi (\LL )\in \Diff_{loc}(P)$
is defined for any $p \in P$ such that $\pi(p)$ is in the 
support of $\phi \in U_{loc}(\gm)$  by 
\begin{equation}
\label{eq:action}
\phi_{\LL}  (p)\equal  {\LL} \cdot p ,
\end{equation}   
where $\cdot$ denotes the $\gm$-action. Note that
$\LL \cdot p$ is uniquely determined  by assumptions.

\begin{prop}
\label{prop:existeuneconnection}
Assume that  $P\to \gm_0$ is  a  principal $G$-bundle over a Lie groupoid
$\gm\toto\gm_0$ and 
$\theta \in \Omega^1 (P)\otimes \Gg$ is
a pseudo-connection. The following statements  are equivalent.
\begin{enumerate}
\item  $\theta  $ is a connection;
\item The $\Gg$-valued one-form $\theta$ is
 basic  with respect to  the $ U_{loc}(\Gamma) $-action;
\item  For each $p\in P$, we have the inclusion
$\hat{A}_{\pi (p)}\subset H_p$, and moreover 
the distribution
$\{H_p\subset T_p P, \ \forall p\in P\}$ is 
preserved under the action of $ U_{loc}(\Gamma) $. Here
$H_p\subset T_p P$ is the horizontal subspace defined by the connection 1-form
 $\theta $, and $\hat{A}_{\pi (p)} $ denotes the subspace of $T_p P$
generated by
the infinitesimal action of the Lie algebroid $A\to \gm_0$.
\end{enumerate}
\end{prop}

The following technical lemma and
Corollary \ref{cor:partialbasic'} will be useful. 

\begin{lem}\label{lem:partialbasic}
Let  $\gm \toto \gm_0$ be a  Lie groupoid 
 and $\omega \in \Omega^k(\gm_0)$  a $k$-form
on $\gm_0$.
Then the following statements are equivalent
\begin{enumerate}
\item $\partial \omega \equal  0$;
\item $\omega$ is basic  with respect to the   $ U_{loc}(\gm )$-action.
\end{enumerate}
\end{lem}
\begin{pf} Note that the Lie algebra of $U_{loc}(\Gamma)$
is the Lie algebra $\Gamma_{loc}(A)$ of local sections of $A$.

1)  $ \Rightarrow$ 2). 
 $\forall r\in \gm$, let $\LL$ be  a  local bisection through
$r$ and  $m\equal \so(r) $.
For all $  \delta_m^1,\dots,\delta_m^k \in T_m \gm_0$, there are  unique  
$ \delta_r^1,\dots,\delta_r^k \in T_r \gm$ tangent to the submanifold $\LL$
 satisfying $\so_*(\delta_r^i)\equal \delta_m^i $
for all $i \in \{1,\dots,k\}$.
It is simple to see that
$\phi_{\LL_*} \delta_r^i \equal \ta_* \delta_m^i  $ for all 
$i \in \{1,\dots,k\}$. Since  $ \partial \omega \equal 0$,  we have
\begin{equation} 
\label{eq:meaningdtheta}
   0 \equal  (\partial \omega ) (\delta_r^1,\dots, \delta_r^k) 
\equal \omega (\delta_m^1,\dots,\delta_m^k) 
-\omega  ( \phi_{\LL*}\delta_m^1,\dots, \phi_{\LL*}\delta_m^k).
 \end{equation}
   Therefore $\phi_{\LL}^* \omega \equal  \omega$.

 $\forall  X \in \Gamma (A)$,  let $\Vec{X}$ be its corresponding right 
invariant vector field  on $\gm $, and $\hat{X}$ the vector field
on $\gm_0$ generated by the infinitesimal action   of $X$ (this is also
denoted by $a(X)$, where $a: A\to T\gm_0$ is the anchor map).
It is known that these vector fields are 
related by the relation $\ta_* \Vec{X}(r)  \equal  \hat{X}(m) $.
  Thus it follows that
 \begin{equation}
\label{eq:VecX}
 \Vec{X} \per \partial \omega
\equal 
\so^* \big( (\so_* \Vec{X}(r))  \per \omega \big)
 - \ta^* \big( (\ta_* \Vec{X}(r))  \per \omega \big) 
\equal - \ta^* \big( (\ta_*\hat{X}(r)  )\per \omega \big)   ,
\end{equation}   
 since $\so_* \Vec{X}(r) \equal 0$.
Since $\ta$ is a surjective submersion, it follows  that
 \begin{equation}
\label{eq:VecX2}
\hat{X}(m)\per \omega \equal   \ta_* \hat{X}(r)  \per \omega \equal 0.
\end{equation}
By Eq.~(\ref{eq:meaningdtheta})
and (\ref{eq:VecX2}),
$\omega$ is a basic form with respect to the $ U_{loc}(\gm)$-action.

2)  $ \Rightarrow$ 1). 
By working backwards using Eq. (\ref{eq:VecX}), one obtains that
$ \Vec{X} \per \partial \omega \equal  0$. {\em I.e.},
$\delta_r \per \partial \omega \equal  0$
if $\so_* \delta_r \equal 0$. Thus it follows that
 $ (\partial \omega )(\delta_r^1 , \cdots ,  \delta_r^k) \equal 0$
whenever $\so_*\delta_r^1 \wedge \dots \wedge \so_*\delta_r^k\equal 0 $.
Similarly,  we have
$(\partial \omega)(\delta_r^1,  \cdots , \delta_r^k)\equal 0$
whenever $\ta_*\delta_r^1 \wedge \cdots \wedge \ta_*\delta_r^k\equal 0 $.

Assume now that $\so_* \delta_r^1 \wedge \cdots 
\wedge \so_*  \delta_r^k \neq 0 $ and  $\ta_* \delta_r^1 
\wedge \cdots \wedge \ta_*  \delta_r^k \neq 0 $. 
Then there exists a local   bisection $\LL$  through $r$ such that
the vectors $\delta_r^1 \cdots , \delta_r^k $ are all tangent to 
the submanifold $\LL$.
Thus, $\forall i\in \{ 1, \cdots , k\}$,  we have
$\phi_{\LL*}\so_* \delta_r^i \equal 
\ta_* \delta_r^i $. It  thus follows that
  $$(\partial \omega)(\delta_r^1 ,  \cdots , \delta_r^k)\equal 
   \omega(\phi_{\LL_*}\so_* \delta_r^1, 
  \cdots , \phi_{\LL_*}\so_* \delta_r^k)
-\omega(\so_*\delta_r^1,   \cdots , \so_*\delta_r^k )\equal 0.$$
Therefore $\partial \omega \equal 0$. This completes the proof of the lemma.
\end{pf}

\begin{cor}\label{cor:partialbasic'}
Let $X \to \gm_0$ be a $\gm$-space
  and $\omega \in \Omega^k (X)$.
Then the following statements are equivalent
\begin{enumerate}
\item $\partial \omega \equal  0$, where $\partial$ is with respect
to the transformation groupoid
$Q :\equal \Gamma \times_{\so,\Gamma_0,\pi} X \toto X$;
\item $\omega$ is basic  with respect to the   $ U_{loc}(\gm)$-action.
\end{enumerate}
\end{cor}
\begin{pf}
Note that  $\omega \in \Omega^k (X)$
is basic with respect to the action of $U_{loc}(\gm)$
if and only if it is basic with respect to the action of $U_{loc}(Q)$.
Then the conclusion follows immediately from 
 Lemma \ref{lem:partialbasic}.
\end{pf}

Now we are ready to  prove Proposition \ref{prop:existeuneconnection}.
\begin{pf}  1) $ \iff $ 2)   follows from Corollary \ref{cor:partialbasic'}. 
 
2) $ \Rightarrow$ 3) is trivial.

3) $ \Rightarrow$ 1)
First, let us fix some notations.
We denote by $\Psi$ the natural map from (an open subset of)
$U_{loc}(\gm) \times P$
to  $Q$ given by 
  \begin{equation}\label{eq:defpsi}
\Psi(\LL, p)\equal  (\LL  \big( \pi(p) \big),p) \in
\gm \times_{p,\gm_0,\pi} P \ \ \ (\equal  Q ).
\end{equation}

 Let $\delta_q \in T_q Q$ be any tangent vector in $Q$ and let
$q(t)\equal \big( \gamma(t), p(t) \big) \in Q$ 
be a path through $q$ such that $\dot{q}(0)\equal \delta_q $.
Let $ \LL (t)$ be a family of local bisections in $ U_{loc}(\gm)$ 
through $ \gamma(t)\in \gm $.  Thus
  $ \Psi( \LL (t), p(t))\equal q(t)$.
Hence we have $\so_* \delta_q \equal \dot{p}(0)$ and
$\ta_* \delta_q \equal  \frac{d}{dt}|_{t\equal 0} (\LL(t)\cdot p(t)) 
\equal  \phi_{\LL*}\dot{p}(0)+ \frac{d}{dt}|_{t=0}
 (\LL (t) \cdot p)$,
where  
$ \LL \equal \LL(0)$ and $ p(0)\equal p$.
Since $\frac{d}{dt}|_{t=0} (\LL (t) \cdot p) \in \hat{A}_{\pi(\LL\cdot p)}$, by assumption we have
  $$  \delta_q \per \partial \theta \equal  
 \dot{p}(0) \per \theta-\phi_{\LL*}\dot{p}(0)  \per \theta.$$   
%
From the decomposition $ T_pP \equal  V_p\oplus H_p  $, we   write 
  $\dot{p}(0)\equal  \hat{X}(p)+ v_{hor}  $,
where $v_{hor}\in  H_p $ is a horizontal vector and
 $ \hat{X} $ is the vector field on $P$  corresponding
to the infinitesimal action of $ X \in \Gg $.
Therefore
  \begin{equation}\label{eq:nul1} 
\delta_q \per \partial \theta
\equal  v_{hor}\per \theta-\phi_{\LL*} v_{hor}\per \theta
 +\hat{X}(p)  \per \theta- \big( \phi_{\LL*} \hat{X}(p) \big)  \per \theta.
\end{equation}
Since $v_{hor}$ is horizontal and the distribution of horizontal subspaces
is preserved under the $U_{loc}(\gm)$-action, we have
  $ v_{hor} \per \theta \equal  \phi_{\LL*} v_{hor}\per \theta\equal 0$. 
Since the  $U_{loc}(\Gamma)$-action commutes with
the  $G$-action,
 we have $\phi_{\LL*} \hat{X}(p)\equal   \hat{X}(\LL \cdot p)$ and therefore 
  \begin{equation}\label{eq:nul2}
\hat{X}(p)\per \theta-\phi_{\LL*} \hat{X}(p) \per \theta  \equal  X-X \equal 0.
\end{equation}
From  Eq.  (\ref{eq:nul1}-\ref{eq:nul2}), it follows
that  $\partial \theta \equal 0$.
\end{pf}

As a special case,  let us consider an  $H$-equivariant
principal (right) G-bundle over $M$ as in Example \ref{ex:equivariant},
we have the following

\begin{cor}
Let $\gm$ be the  transformation groupoid $H\times M\toto M$, where
$H$ acts on $M$ from the left.
 And let  $P\to M$  be
an $H$-equivariant principal (right) $G$-bundle over $M$
considered as a $G$-bundle over $\gm $. A one-form 
$\theta  \in \Omega^1 (P)\otimes \Gg$ defines a connection
in the sense of Definition \ref{def:connection} 
if and only  if $\theta $ is a connection $1$-form for the
 $G$-bundle $P\to M$ and is basic with respect to the $H$-action.

In particular, if $H$ acts on $M$ freely and properly, a connection
on the principal $G$-bundle $P\to M$ over $\gm$ is
equivalent to a connection on the (ordinary) principal $G$-bundle
$P/H\to M/H$.
\end{cor}
\begin{pf} 
Let $Q$ be the transformation groupoid $\gm \times_M P\toto P$.
It is simple to see that $Q$ is isomorphic to 
$H\times P\toto P$. By $\so$ and $\ta$ we denote
the source and target maps.
$\forall X\in {\mathfrak{h}}$,  by $\Vec{X}$ we denote
its   corresponding right-invariant vector field on $H$ and by
$\hat X$ the vector field on $P$ induced by the infinitesimal action
of $X$. Thus $\forall p\in P, h\in H$ and $v\in T_pP$, we have
$$(\so^*\theta)(\Vec{X}(h),v)\equal \theta(v)\ \mbox{and }
(\ta^*\theta)(\Vec{X}(h),v)\equal \theta(\hat X_{hp})+\theta(\phi_{h*}v),$$
where $\phi_h$ denotes the $H$-action on $P$.
Therefore $\partial\theta\equal 0$ if and only if
$\theta(\phi_{h*}v )\equal \theta(v) , \ \forall h\in H, \  v\in TP$ and
$\theta(\hat X_p)\equal 0, \ \forall X\in \Hh$, {\em i.e.}
$\theta$ is $H$-basic.
\end{pf}

Next we study  some obvious obstructions to the
existence of  a connection. 
Recall that a principal $G$-bundle $P$ over
$\Gamma  \toto \Gamma_0$ induces 
a Lie groupoid homomorphism $\rho: \Gamma  \to \frac{P \times P}{G} $
given by Eq.~(\ref{eq:Q1}). For any $m\in \Gamma_0$,
we denote by  $I_m\equal \gm_m^m$
the isotropy group at $m \in \Gamma_0$ and $(I_m)_0$ its connected component
of the identity.  A groupoid is said to be a {\it foliation groupoid}
if all its isotropy groups are finite groups \cite{Moe}.

\begin{prop}
If the principal $G$-bundle $P$ over
$\Gamma  \toto \Gamma_0$ admits a connection, then
\begin{enumerate}
\item $(I_m)_0$  acts on $P$ trivially,
\item $(I_m)_0 \subseteq \Ker \rho $, and
\item if $ \Gamma$ is proper, the image $ \rho(\gm)$  is an immersed 
foliation subgroupoid of 
the gauge groupoid $\frac{P \times P}{G} $.
\end{enumerate}
\end{prop}

\begin{numrk} We  note that  $\rho(\gm) \toto \gm_0$
is neither  a submanifold of $\frac{P \times P}{G}$  nor
 a Lie groupoid in general. 
\end{numrk}

\begin{pf}
1) For any $p \in P$, the infinitesimal action of the Lie
 algebra ${\mathfrak i}_m $
of $ I_m$, where $m\equal \pi (p)$,
 defines a subspace $ \hat{{\mathfrak i}}_m \subset T_pP$.
 On the  one hand, since $\hat{{\mathfrak i}}_m \subset \hat{A}_m$,  by 
 Proposition \ref{prop:existeuneconnection}(3), we have 
    $\hat{{\mathfrak i}}_m \subset H_p $,
where $ H_p$ is the horizontal tangent space at $p$.
  On the other hand, since the action of  $ I_m$ maps $ \pi^{-1}(m)$
to itself, $\hat{{\mathfrak i}}_m$ must be  a subspace
 of  $ T_p (\pi^{-1}(m) )$ or  $V_p$, the vertical tangent
 space at $p$,  {\em i.e.}, $\hat{{\mathfrak i}}_m \subseteq  V_p  $.
 Since $H_p  \cap V_p\equal 0$, it thus follows that
$\hat{{\mathfrak i}}_m\equal 0 $.
Hence, the Lie group $ (I_m)_0$ acts on $P$  trivially.

2) is just a rephrase of 1).

3) Since $\rho$ is the identity map when being restricted
to the unit space, 
the isotropy group $J_m $ of $\rho(\gm)$ over $m \in \gm_0$ is given by
   $  J_m \equal  \rho (I_m)$.
 Since $(I_m)_0 \subset \ker (\rho) $,  it thus follows that
$ \rho$ passes to the 
quotient $J_m \equal  {\rho} \big(\frac{I_m}{(I_m)_0}\big)$.
    The groupoid $\gm$ is  proper, so $ I_m$ must be compact 
and therefore $\frac{I_m}{(I_m)_0} $ is discrete.
It thus follows that
the group $J_m$ is a finite group.
Hence the groupoid   $ \rho(\gm)$ is a
foliation groupoid.
\end{pf}

\begin{numex}
\label{ex:3.10}
Let $G$ be a Lie group. Then
$G\to \cdot$ can be considered as a principal
$G$-bundle over the groupoid $G\toto \cdot$, where the groupoid 
$G\toto \cdot$ acts on $G$ by left translation. It is clear that
connection does not exist unless $G$ is discrete.
\end{numex} 

\subsection{Pull-back connections}

In this subsection, we show that there is a natural notion of
pull-back connections under generalized homomorphisms which
generalizes the usual pull-back connections of principal
bundles over manifolds.

Let $P$ be a principal  $G$-bundle over $\Gamma \toto \Gamma_0 $,
and  $ \phi:\Gamma'_0 \stackrel{\tau }{\leftarrow}
 Z \stackrel{\sigma}{\rightarrow} \Gamma_0$
a generalized homomorphism from $\Gamma' \toto \Gamma'_0$
to $\Gamma \toto \Gamma_0$. By  $P'$, we denote the pull-back
principal $G$-bundle over $\Gamma' \toto \Gamma'_0$ constructed
as in Proposition \ref{prop:pullback}. 

We need to introduce some notations.
Let $\tsigma: Z\times_{\Gamma_0}P \to P' \ \ \big(\equal
(Z\times_{\Gamma_0}P)/\gm \big)$
 be the natural projection and $\ttau: Z\times_{\Gamma_0}P\to P$
  the projection map
to the second component.

\begin{prop}\label{prop:connecetgeneral}
If $\theta \in   \Omega^1 (P)\otimes \Gg$  is a connection 
for the principal $G$-bundle $P$ over $\Gamma \toto \Gamma_0 $, then
\begin{enumerate}
\item there is a unique one-form $\theta' \in   \Omega^1 (P')\otimes \Gg$ 
satisfying the condition:
\begin{equation}
\label{eq:pratique}
\tsigma^* \theta' \equal \ttau^* \theta; 
\end{equation} 
\item $\theta' $ defines a connection on the pull-back principal $G$-bundle
 $P'$ over $\Gamma' \toto \Gamma'_0$, which is called the pull-back
connection and is denoted by $\phi^* \theta$;
\item  the curvature $\Omega' $ of $\theta'$
 and the curvature $\Omega$   of $\theta$ are related by
 \begin{equation} \label{eq:omegaestnul} \tsigma^* \Omega'
\equal \ttau^* \Omega ;
\end{equation}
\item $\phi^* \theta $ is flat if $\theta $ is flat;
\item  if  $\psi$ and $\phi$ are generalized homomorphisms from
$\gm''\toto \gm''_0 $ to  $\Gamma' \toto \Gamma_0' $ and
from $\Gamma' \toto \Gamma_0' $ to $\gm \toto \gm_0$ respectively, then
$$(\phi \smalcirc \psi)^*\theta \equal \psi^* (\phi^* \theta ).$$
\end{enumerate}
\end{prop}  

\begin{pf}  1) Note that $ Z \times_{\gm_0} P \to \gm_0$ is a $ \gm$-space
and hence admits a $U_{loc}(\Gamma)$-action.  Moreover we have 
\begin{equation} 
\label{eq:isom} 
(Z \times_{\gm_0} P) /U_{loc}(\gm) \simeq  (Z \times_{\gm_0} P) /\gm
\equal P'   .
\end{equation}
Let $ \theta_Z\equal  \ttau^* \theta \in \Omega^1(Z \times_{\gm_0} P) 
\otimes \Gg$.
Since the projection $\ttau: Z \times_{\gm_0} P \to P $  commutes with
the $ U_{loc}(\gm)$-action, and  $ \theta$ is basic with respect to
the $ U_{loc}(\gm)$, the  1-form $\theta_Z$ is also basic. Hence there exists 
an unique $\theta' \in \Omega^1(  P') \otimes \Gg$  such that
$ \tsigma^* \theta' \equal \theta_Z$.

2) The triple $(Z\times_{\Gamma_0}P,\ttau,\tsigma)$
 defines a generalized homomorphism from $Q' \toto P' $ to
$Q \toto P$, where $Q' \toto P' $ is
the transformation groupoid of the principal $G$-bundle $P'\to \gm_0'$ 
over $\gm' \toto \gm'_0$.
Let   $ Q'[Z\times_{\Gamma_0} P]~\toto~Z\times_{\Gamma_0}~P $ be
the pull-back groupoid
of $ Q' \toto P'$ via $ \tsigma :~Z\times_{\Gamma_0}~P~\to~P'$.
By abuse of notation, we denote by $ \tsigma$  and $\ttau$  the homomorphisms
from $ Q'[Z\times_{\Gamma_0} P] \toto Z\times_{\Gamma_0} P $
to $ Q' \toto P'$ and from
$ Q'[Z\times_{\Gamma_0} P] \toto Z\times_{\Gamma_0} P $
to $ Q \toto P$,   respectively, defined as in Eqs.~(\ref{eq:deftau})
and (\ref{eq:defsigma}):
$$ \begin{array}{ccccc} 
              Q' & \stackrel{\tsigma}{\leftarrow}                   
& Q'[Z\times_{\Gamma_0} P]           &\stackrel{\ttau}{\rightarrow}     & Q \\
 \downdownarrows &
& \downdownarrows & & \downdownarrows \\
              P' &  \stackrel{\tsigma}{\leftarrow}
& Z \times_{\gm_0} P         &\stackrel{\ttau}{\rightarrow} &P\\
       \end{array} $$
We can now prove 2) easily. Since $\tsigma$ and $\ttau$ are strict
groupoid  homomorphisms,
we have,  
\begin{equation} \label{eq:re1} \tsigma ^*\partial \theta'
\equal  \partial \tsigma^* \theta'\equal  \partial \theta_Z,
 \end{equation}
 and on the other hand,  
 \begin{equation} \label{eq:re2} 0 \equal  \ttau^* \partial \theta \equal  
  \partial \theta_Z  .\end{equation}
From Eqs.  (\ref{eq:re1}-\ref{eq:re2}), it follows that
  $ \tsigma^* \partial \theta' \equal  0$.
Since $  \tsigma : Q'[Z\times_{\Gamma_0} P] \to Q'$
is a surjective submersion, we must have $ \partial \theta '\equal  0$.
This proves 2).

3) follows from Eq. (\ref{eq:pratique}).  4) is straightforward.
5) is a direct verification and is left to the reader.
\end{pf}



\begin{numrk} 
If $ \phi:\Gamma' \to \gm$ is a groupoid strict homomorphism,
then the pull-back can be described more explicitly.
In this case $P'\equal \phi^* P  \equal  \gm' \times_{\phi,\gm_0,\pi} P  $ 
and $ \theta'\equal  \pr^* \theta$,
 where $ \pr: \gm' \times_{\phi,\gm_0,\pi} P \to P$
 is the projection on the second component.
\end{numrk}

\begin{cor} 
\label{cor:3.15}
If  $\gm'$ and  $\gm$ are Morita equivalent groupoids,
there is an equivalence of categories of principal
$G$-bundles with connections over $\gm'$ and  $\gm$.
There is also an equivalence of categories of principal
$G$-bundles with flat connections over $\gm'$ and  $\gm$.
\end{cor}

Therefore, this allows us to speak about
 {\em connections} and {\em flat connections}
of principal $G$-bundles over a differential stack.
As we see from Example \ref{ex:3.10}, connections may
not always exist. However, for orbifolds,
we will show that they always exist.

\begin{them}
\label{prop:orbifoldconnection}
Any principal $G$-bundle over an orbifold  admits a connection.
\end{them}
\begin{pf} It is well-known \cite{Moe} that  an 
 orbifold can be represented by a  proper \'etale groupoid
$\gm \toto \gm_0$.

Let $\theta \in \Omega^1(P) \otimes \Gg$ be a pseudo-connection
of the principal $G$-bundle $P \stackrel{\pi}{\to} \gm_0$.
For any $\gamma\in \Gamma$, there is always a unique
local bisection through $\gamma$ since $\Gamma$
is an \'etale groupoid.
We denote by   $\phi_{\gamma}$ the local diffeomorphism on
$P$ induced  by this local bisection. Note that
 $\phi_{\gamma}$ is defined in an open neighborhood of
any $p\in \pi^{-1}(\so(\gamma))$.

Since $\Gamma\toto \gm_0$ is a Lie groupoid, there exists
 a (right) Haar system, denoted by $\lambda\equal (\lambda_x)_{x\in \gm_0}$,
where $\lambda_x$ is a measure with support on $\Gamma_x\equal \so^{-1}(x)$ such
that
\begin{itemize}
\item for all  $f\in C_c^\infty(\Gamma)$,
$x\mapsto \int_{\gamma \in\Gamma_x} f(\gamma)\,\lambda_x(d\gamma)$ is
smooth, and
\item the right translation  by $R_{\gamma'}: \gm_y \to \gm_x$ (where
$y\equal \ta(\gamma')$
and $x \equal  \so(\gamma')$ for all $\gamma' \in \Gamma$) preserves the measure, {\em i.e.}
\begin{equation}\label{eq:haar}  \int _{\gamma \in\Gamma_y} f(\gamma
\gamma')\,\lambda_y(d\gamma)\equal 
 \int _{\gamma \in\Gamma_x} f(\gamma )\,\lambda_x(d\gamma).
 \end{equation}
\end{itemize}

Let us recall that a smooth function $c\colon \gm_0 \to \rr_+$ is
called a {\em  cutoff} function if
\begin{itemize}
\item[(i)] for any $x\in \gm_0$, $\int_{\gamma\in \Gamma_x}
c(\ta(\gamma))\,\lambda_x(d\gamma)\equal 1$;
\item[(ii)] for any $K\subset \gm_0$ compact, the support of
$(c\smalcirc \ta)_{|\Gamma_K}$ is compact.
\end{itemize}
It is known \cite[Proposition~6.7]{tu99}  that  a cutoff  function
exists if and only if $\Gamma \toto
\gm_0$
is proper.
We can now define a connection $\tilde{\theta}$ on $P \stackrel{\pi}{\to}
\gm_0 $
by

   \begin{equation}\label{eq:applatie}  \tilde{\theta}_p\equal 
\int_{\gamma \in\Gamma_{\pi(p)}}
c(\ta(\gamma))
\, \, \phi_{\gamma}^* \theta_{\gamma \cdot p}\, \,
\,\lambda_{\pi(p)}(d\gamma)
  .\end{equation}

Since $c(\ta(\gamma))$  has a compact support in $\Gamma_{\pi(p)}$,
the integral in  Eq.~(\ref{eq:applatie}) is well-defined.
Let us check that $\tilde{\theta}$ is indeed a connection.
It is easy   to see,   by  Eq.~(\ref{eq:applatie}), 
that  $\tilde{\theta}$ is $G$-invariant. For any $ X \in \Gg$,

\begin{equation}\label{eq:cond2'}  
\hat{X}(p) \per \tilde{\theta}_p\equal  \int_{\gamma \in\Gamma_{\pi(p)}}
c(\ta(\gamma))
\, \, \big( \phi_{\gamma*}\hat{X}(p)\big) \per \theta_{\gamma \cdot p}\, \,
\,\lambda_{\pi(p)}(d\gamma).
  \end{equation}

Since $\phi_{\gamma*}\hat{X}(p)\equal \hat{X}(\gamma \cdot p) $
and $\hat{X}(\gamma \cdot p) \per \theta_{\gamma \cdot p}\equal   X $,
Eq.~(\ref{eq:cond2'}) implies that $\hat{X}(p)\per \tilde{\theta_p}
\equal  X$.

We check  that $\tilde{\theta}$ is basic with respect to
the  $U_{loc}(\gm)$-action.
Since $\gm$ is \'etale, we only need to check
that $\tilde{\theta}$ is invariant under the  $ U_{loc}(\gm)$-action, or
 $\phi_{\gamma'}^* \tilde{\theta}\equal \tilde{\theta}$ for any
$\gamma' \in \gm$.

By Eq.~(\ref{eq:applatie}), we have, for any $\gamma'$ with $\ta(\gamma')
\equal  \pi(p)$,
  \begin{equation}\label{eq:cond7} \phi_{\gamma'}^* \tilde{\theta}_p
\equal    \int_{\gamma \in\Gamma_{\pi(p)}} \,
c(\ta(\gamma)) \,\,\,
\phi_{\gamma'}^*\phi_{\gamma}^* \theta_{\gamma \cdot p} \,\,\,
\,\lambda_{\pi(p)}(d\gamma) .\end{equation}
Now, using the relations
$\phi_{\gamma'}^* \smalcirc \phi_{\gamma}^*   \equal  \phi_{\gamma  \gamma'}^*
$, $ c(\ta(\gamma))  \equal  c(\ta(\gamma  \gamma'))$,
 $\theta_{\gamma \cdot p} \equal   \theta_{\gamma \gamma' \cdot p'}$,
where $p'\equal \gamma'^{-1} \cdot p$, we have \begin{equation}\label{eq:cond8} 
 \phi_{\gamma'}^* \tilde{\theta}_p
\equal    \int_{\gamma \in\Gamma_{\pi(p)}}
c(\ta(\gamma  \gamma'))\, \,
\phi_{\gamma  \gamma'}^*  \theta_{\gamma  \gamma' \cdot p'}
\, \, \lambda_{\pi(p)}(d\gamma).
\end{equation}
Since the Haar measure  $\lambda_{\pi(p)}(d\gamma)$ is invariant under
the right translation, Eqs.~(\ref{eq:cond8}) and 
(\ref{eq:haar}) imply that
\begin{equation}
\label{eq:cond9}  
\phi_{\gamma'}^*  \tilde{\theta}_p
\equal    \int_{\gamma  \in\Gamma_{\pi(p')}} \,\,
c(\ta(\gamma ))\,\,
\phi(\gamma)^* \theta_{\gamma  \cdot p'}\,\, \lambda_{\pi(p')}(d\gamma) .
\end{equation}
The right hand side of Eq.~(\ref{eq:cond9}) is $\tilde{\theta}_{p'}$
by definition. Therefore we have proved that 
$\tilde{\theta}_{p'}\equal  \phi_{\gamma}^* \tilde{\theta}_p$.
 By Corollary \ref{cor:partialbasic'}, we have
$\partial \tilde{\theta} \equal  0$. This completes the proof.
\end{pf}

\subsection{Holonomy map}

As usual, a flat connection on a principal $G$-bundle over
a groupoid $\gm $  induces a representation of 
the fundamental group of $\gm $ to $G$ via the holonomy map.

Let us recall the definition of  the fundamental group of
a topological groupoid.
We need a few preliminaries.

\begin{prop}
\label{pro:point}
Let $\Gamma$ be a topological groupoid. There is a bijection between
\begin{itemize}
\item[(a)] equivalence classes of 
generalized homomorphisms from $\{pt\}$ to $\Gamma$, and
\item[(b)] elements in the quotient space $\Gamma_0/\Gamma$.
\end{itemize}
\end{prop}

\begin{pf}
A generalized homomorphism from $\{pt\}$ to $\Gamma$ is given by a space
$Z$, a continuous map $p\colon Z\to\Gamma_0$ and a right action of $\Gamma$
on $Z$ whose momentum  map is $p$, such that $Z/\Gamma\simeq \{pt\}$.
The element of $\Gamma_0/\Gamma$ corresponding to that generalized
homomorphism is the class of $p(z)$ for any $z\in Z$.

Conversely, if $\bar x\in \Gamma_0/\Gamma$, choose a representative
$x\in \Gamma_0$. Then $Z\equal  \alpha^{-1}(x) $ 
together with $p = \so:Z\to\Gamma_0$ and the multiplication
on the right as the right action on $Z$,
is a generalized homomorphism from
$\{pt\}$ to $\Gamma$.
It is clear that if $x, \ y\in \Gamma_0$ are in the
  same groupoid orbit their corresponding
generalized homomorphisms are   equivalent.
\end{pf}

We call a generalized homomorphism  
$\{pt\} \stackrel{}{\leftarrow}X\stackrel{}{\to} \Gamma_0$
a {\em point} in $\Gamma$. A 
{\em pointed groupoid} is a pair $(\Gamma, X)$,  where
$\Gamma$ is a groupoid, and $\{pt\} \stackrel{}{\leftarrow}X\stackrel{}{\to} \Gamma_0$
 a point in $\Gamma$. If $\Gamma'_0
\stackrel{\rho}{\leftarrow}Z\stackrel{\sigma}{\to} \Gamma_0$ is a generalized
homomorphism from a groupoid $\Gamma'$ to a groupoid $\Gamma$,
and $X'$ is a point in $\Gamma'$, then one can define
the image $Z\smalcirc X'$, a point in $\Gamma$ by using the composition of
generalized homomorphisms
$$ \{pt\} \stackrel{X'}{\to}\Gamma' \stackrel{Z}{\to} \Gamma.$$

Suppose that $(\Gamma,X)$ and $(\Gamma',X')$ are two pointed groupoids. Then
a pointed generalized homomorphism from $(\Gamma',X')$ to $(\Gamma,X)$ is a
pair $(Z,\varphi)$, where $Z$ is a generalized homomorphism
from $\Gamma'$ to $\Gamma$ and $\varphi$ is an equivalence from
$Z \smalcirc X'$ to $X$.
It is clear that pointed groupoids
with  pointed generalized homomorphisms form a category.

%

\begin{rmk}
Consider the special case when $\Gamma'$ is
a topological space $M$, and $X'$ is a basepoint $ m_0\in M$.
For any $x_0\in\Gamma_0$, let 
$\{pt\} \stackrel{}{\leftarrow}\ta^{-1}(x_0)\stackrel{}{\to} \Gamma_0$
be the point  in $\gm$ as in the proof of Proposition \ref{pro:point},
which  is denoted by $x_0$ by abuse of notation.   Then
 a pointed generalized homomorphism from $(M,m_0)$ to $(\Gamma, x_0)$
is equivalent to giving a generalized homomorphism
$M\stackrel{\tau}{\leftarrow} Z\stackrel{\sigma}{\to} \Gamma_0$
and a point $z_0\in Z$ such that $\tau(z_0)\equal m_0$
and $\sigma(z_0)\equal x_0$. Here $z_0\in Z$
is the unique element corresponding to the inverse image of
$\epsilon(x_0)$ under the isomorphism $\phi: \tau^{-1}(m_0 )\to \ta^{-1} (x_0)$
(where $\epsilon : \gm_0 \to \gm$ is the unit map).
\end{rmk}

We say that two pointed generalized homomorphisms $(Z_0,\varphi_0)$
and $(Z_1,\varphi_1)$ from $(M,m_0)$ to $(\Gamma,X)$
are {\em homotopic} if there exists a pointed
generalized homomorphism $(Z,\varphi)$ from $((M\times [0,1])/(\{m_0\}
\times [0,1]), m_0 )$ to $(\Gamma,X)$ such that
\begin{itemize}
\item there exist equivalences $i_0\colon Z_0\to Z_{|M\times\{0\}}$
and $i_1\colon Z_1\to Z_{|M\times\{1\}}$;
\item via the identifications $i_0$ and $i_1$ above, the restriction
of $\varphi$ to $Z_0\smalcirc X'$ and $Z_1\smalcirc X'$ is equal to
$\varphi_0$ and $\varphi_1$, respectively.
\end{itemize}

We denote by $[(M,m_0),(\Gamma,X)]$ the set of homotopy classes
of pointed generalized homomorphisms from $(M,m_0)$ to $(\Gamma,X)$.

\begin{rmk}
In the non-pointed situation, the set of homotopy classes
of generalized homomorphisms from a space $M$ to a groupoid $\Gamma$
is also called the set of concordance classes of principal
$\Gamma$-bundles over $M$, or $\Gamma$-structures on $M$.
\end{rmk}

Let $(M',m'_0)$ be another topological space. By
$(M,m_0)\vee (M',m'_0)$, we denote the pointed space $(M\cup M')/(m_0\sim m'_0)$
with the basepoint $m_0\equal m'_0$. It is clear that if $(Z,\varphi)$
(resp. $(Z',\varphi')$) is a pointed generalized homomorphism from
$(M,m_0)$ (resp. $(M',m'_0)$) to $(\Gamma,X)$ then one can form
a pointed generalized homomorphism from $(M,m_0)\vee (M',m'_0)$ to
$(\Gamma,X)$. For $M\equal M'\equal S^n$, using the usual map
$(M,\ast)\to (M,\ast)\vee (M,\ast)$ we get a group structure on
$\pi_n(\Gamma_{\com} ,X):\equal  [(S^n,\ast),(\Gamma,X)]$.
It is clear that if $X$ and $Y$ are equivalent 
points in $\Gamma$, then $\pi_n(\Gamma_{\com} ,X)$ is
isomorphic to $\pi_n(\Gamma_{\com} ,Y)$.
 The group $\pi_1(\Gamma_{\com} ,X)$ is 
called the fundamental group of $\Gamma$.
For any $x\in \gm_0$, consider the point 
$\{pt\} \stackrel{}{\leftarrow}\ta^{-1}(x)\stackrel{}{\to} \Gamma_0$
in $\gm$ as in the proof of Proposition \ref{pro:point}. We  denote
by $\pi_1(\Gamma_\com, x)$  its corresponding fundamental group.
Since the isomorphism class of $\pi_1(\Gamma_\com, x)$
only depends on the class $\bar{x}\in \Gamma_0/\Gamma$,
we also denote it by  $\pi_1(\Gamma, \bar{x})$ as well.

\begin{rmk}
We recently learned that the fundamental group  of a stack is also
being introduced by Behrang Noohi \cite{Noo}
\end{rmk}

It will become clear below (Proposition~\ref{prop:BGamma})
that if $x$ and $y$ are in the same
path-connected component of $\Gamma_0$ then $\pi_1(\Gamma_\com,x)$
and $\pi_1(\Gamma_\com,y)$ are isomorphic (but the isomorphism may be not canonical).
As a consequence, we have 

\begin{prop}
Let $\Gamma$ be a topological groupoid.
Suppose that $\Gamma_0/\Gamma$ is path connected,
 then $\pi_1(\Gamma_\com , x)$
is independent of the choice of $x$ up to isomorphism.
\end{prop}

\begin{rmk}
It is clear by construction that if $(\Gamma,X)$ and
$(\Gamma',X')$ are two Morita-equivalent pointed groupoids
then $\pi_n(\Gamma_\com,X)$ and $\pi_n(\Gamma'_\com,X')$ are
isomorphic.
\end{rmk}

\begin{rmk}
By a generalized path on $\Gamma$, we mean a generalized homomorphism
from an interval $[a,b]$ to $\Gamma$. If $P$ and $P'$ are generalized
paths from $[a,b]$ and $[b,c]$ to $\Gamma$ respectively,
and if an equivalence $\psi$ from $P\smalcirc  b$ to $P'\smalcirc b$ is given, then
the composition of $P$ and $P'$ is the generalized path
$P''$ from $[a,c]$ to $\Gamma$ defined by $(P\cup P')/(p\sim \psi(p))$.
\end{rmk}

\par\bigskip
We will see below that the homotopy groups of $\Gamma$ are in fact
isomorphic to the homotopy groups of the classifying space
$B\Gamma$, which is the fat geometric 
realization of the simplicial space $\gm_\com$.


%

Recall that the space $B\Gamma$ is defined by $\underrightarrow{\lim_n} (B\Gamma)_n$,
 where
$$(B\Gamma)_n\equal \left( \coprod_{0\le k\le n} \Gamma_{k} \times \Delta_k\right) / \sim, $$
and  the equivalence relation $\sim$ is generated
by 
$$\big(\gamma, \tilde{\epsilon}^k_i(t) \big)
\sim \big(\epsilon_i^k (\gamma),t \big) 
\ \ \ {\mbox{and}} \ \ \ \big( \gamma ,\tilde{\eta}_i^k (t') \big)
\sim \big(\eta_i^k (\gamma), t' \big)  \ \ \ \forall \gamma \in \gm_k, 
\forall t\in \Delta_{k-1},
 \forall  t' \in \Delta_{k+1}.$$ 
Here ${\epsilon}^k_i$
and ${\eta}_i^k $ are the face and degeneracy maps of
$\gm_\com$ defined by Eqs. (\ref{eq:face}-\ref{eq:deg}), while
$\tilde{\epsilon}^k_i$ and $\tilde{\eta}_i^k $ are the face and degeneracy 
maps of $\Delta_\com$ defined by Eqs. (\ref{eq:face1}-\ref{eq:deg1}).
It is known that   $B\gm$ is naturally endowed with
a topology   \cite{Dupont,Segal}.

\begin{prop}\label{prop:BGamma}
Let $(M,m_0)$ be a locally compact $\sigma$-compact
space, $\Gamma$ a topological
groupoid and $x_0\in \Gamma_0$.
Then there is a canonical isomorphism
$$[(M,m_0),(\Gamma, x_0)] \cong [(M,m_0),(B\Gamma,x_0)].$$
\end{prop}

(In the theorem above, $\Gamma_0$ is identified to a subspace of $B\Gamma$).

\begin{pf}
The proof is almost identical to that of \cite[Theorem I.7]{hae70} or
\cite{bl70}, where the proposition is proved in the unpointed case.
\end{pf}

\begin{cor}
For all $n\in \nn$ and $x\in \gm_0$,
 $\pi_n(\Gamma, x)\cong \pi_n(B\Gamma,x)$.
\end{cor}

\par\bigskip
Before we introduce the holonomy map, we need to characterize
geometrically flat $G$-principal bundles over 
a groupoid $\Gamma\toto\Gamma_0$
({\it i.e.} $G$-principal bundles which admit a flat connection).

\begin{prop}\label{prop:flat}
Let $P$ be a $G$-principal bundle over $\Gamma\toto \Gamma_0$.
The following are equivalent:
\begin{itemize}
\item[(i)] $P$ is flat;
\item[(ii)] there exist an open cover $(U_i)$ of $\Gamma_0$ and a
local trivialization $\phi_i\colon U_i\times G \tilde{\to} P|_{U_i}$
  such that the transition function 
$\psi_{ij}\colon \Gamma^{U_i}_{U_j}\to G$ defined by the equation
\begin{equation}
\label{eq:tran}
\gamma\cdot\phi_j(\so(\gamma),g) \equal  \phi_i(\ta(\gamma),\psi_{ij}(\gamma)g)
\end{equation}
is a constant map for any $\Gamma^{U_i}_{U_j}$;
\item[(iii)]  As a  $G$-principal bundle over $\Gamma$, $P$
  is induced from a $G^d$-principal bundle $P'$ under the
natural group homomorphism $G^d \to G$, 
 where $G^d$ is the group $G$ endowed with the discrete topology.
\end{itemize}
\end{prop}

\begin{pf}
(i) $\implies$ (ii): Since $P$ is flat as a $G$-principal
bundle over $\Gamma_0$, there exists a local trivialization $\phi_i
\colon U_i\times G\to P|_{U_i}$ such that $\phi_i^*\theta$
is the trivial connection on $U_i\times G$ for all $i$
({\it i.e.} $\phi_i(U_i\times \{g\})$ is horizontal for all $g\in G$).
Since $\del\theta \equal 0$, horizontal sections
are invariant under the $U(\gm )$-action according
to Lemma \ref{lem:partialbasic}. It thus follows 
that $\psi_{ij} (\gamma )$ must be independent of
$\gamma$. 

\par\medskip

(ii) $\implies$ (iii):
Let $P'$ be the $G^d$-principal bundle over $\Gamma$
defined by the cocycle $\{\psi_{ij}\}$. More explicitly,
let $\gm [ U_i]\toto \coprod U_i$ be the pull back
groupoid of $\gm \toto \gm_0$ via the open covering
$\coprod U_i\to \gm_0$. Then $\gm [ U_i] \equal \coprod 
\Gamma^{U_i}_{U_j}$ and $\coprod \psi_{ij}:  \coprod\Gamma^{U_i}_{U_j}
\to G$ is a  groupoid homomorphism, which is constant
on each $\Gamma^{U_i}_{U_j}$. Therefore, it defines a 
continuous groupoid homomorphism from   $\gm [ U_i]$ to $G^d$.
Composing the Morita equivalence from $\gm$ to $\gm [ U_i]$
with this homomorphism, one obtains a generalized 
homomorphism from $\gm$ to $G^d$, which corresponds 
to a $G^d$-principal bundle $P'$ over $\Gamma$.

\par\medskip  

(iii) $\implies$ (i) By assumption, there is a 
local trivialization $\phi_i\colon U_i\times G\to P|_{U_i}$
  such that the transition functions
$\psi_{ij}\colon \Gamma^{U_i}_{U_j}\to G$ defined by Eq.  (\ref{eq:tran})
 are constant on each $\Gamma^{U_i}_{U_j}$.  Therefore 
$\phi_i(U_i\times \{g\})$, for all $g\in G$, induces a
well-defined   horizontal  distribution on $P$, which 
is clearly a flat connection according to
Proposition  \ref{prop:existeuneconnection}.
\end{pf}

\par\bigskip
We are now ready to prove the main theorem of this subsection:
\begin{them}
Let $P\stackrel{\pi}{\to} \Gamma_0$ be a flat $G$-principal bundle over a Lie groupoid
$\Gamma\toto\Gamma_0$, $x\in\Gamma_0$ and $p\in P$ such that $\pi(p)=x$. There is a group
homomorphism
$${\mathrm{Hol}}_{\Gamma,x,p}\colon\pi_1(\Gamma_\com,x)\to G,$$
called the holonomy map, defined as follows.
For every generalized pointed loop $(Z,z_0)$ from
$(S^1,*)$ to $(\Gamma, x)$,
let  $P_{S^1}\equal Z\times_\Gamma P'$  be the
pull back  $G^d$-principal bundle over $S^1$, where $P'$ 
is the corresponding $G^d$-principal bundle over $\gm$ as in Proposition
\ref{prop:flat}.  Let $p'\equal (z_0,p)$.
Then
$${\mathrm{Hol}}_{\Gamma,x,p}(Z,z_0)
:\equal {\mathrm{Hol}}_{S^1,*,p'}(\mbox{Id}_{S^1}).$$
In other words, if $f(t)$ is the horizontal lift (with respect to the pull back connection)
 of the path $\mbox{Id}_{S^1}$ on $P_{S^1}$ starting at $p'$, then
${\mathrm{Hol}}_{\Gamma,x,p}(Z,z_0)$ is the element $g\in G$ satisfying
$f(1)\equal f(0)g$.
\end{them}

\begin{pf}
It remains to prove that the holonomy map
only depends on the homotopy class of the generalized path.
This clearly follows from Proposition~\ref{prop:flat}(iii)
and  the continuity of the holonomy map.
\end{pf}

\begin{rmk}
As usual, if $x$ and $y$ are in the same path-connected
component, then $\pi_1(\Gamma_\com,x)$ and $\pi_1(\Gamma_\com,y)$
are isomorphic, the isomorphism being well-defined up to
an inner automorphism. Via this identification, there exists $g\in G$
such that ${\mathrm{Hol}}_{\Gamma,y}\equal \mbox{Ad}_g\smalcirc
{\mathrm{Hol}}_{\Gamma,x}$.
\end{rmk}

\section{Chern-Weil  map for $G$-differential simplicial algebras}

In this section,  we present an algebraic construction
of Chern-Weil  map for $G$-differential simplicial algebras.
By an algebra, we always mean an ${\mathbb N}$-graded
super-commutative algebra.

\subsection{Simplicial algebras}

We start with recalling the definition of 
$G$-differential algebras \cite{guillemin}.

\begin{defn}
 A $G$-differential algebra is an ${\mathbb N}$-graded differential algebra
$(A, d)$ ($d$ is of degree $1$)  equipped with a $G$-action and
a linear map $\ii $, called the contraction, from
 $\Gg$ to the space of derivations of degree $-1$    satisfying 
\begin{itemize}
\item  the action of $G$ preserves the ${\mathbb N}$-graded differential 
algebra structure;
\item the following identities hold:
 $ \forall X, Y\in \Gg$ and $\forall g\in G$
\begin{eqnarray}
&&\ii_X \ii_Y+ \ii_Y\ii_X\equal 0, \\
&&g\smalcirc \ii_X \smalcirc g^{-1}\equal \ii_{Ad_g (X)},\\
&&\ii_{X}d+d\ii_{X}\equal  L_{X},
\end{eqnarray}
where $L_{X}$ denotes the infinitesimal $\Gg$-action.
\end{itemize}
An element $a\in A$ is called {\em basic} if $L_X  a\equal 0$ and
$\ii_X a\equal 0$ for
any $X\in \Gg$.
\end{defn}

A standard example of $G$-differential algebra is the Weil algebra.

\begin{numex}
({\bf Weil algebra}) 
Let $\Gg$ be the Lie algebra of a Lie group
$G$, and let  $W(\Gg)\equal S(\Gg^*) \otimes \land \Gg^* $ be the tensor
product of $S(\Gg^*)$ and $\land \Gg^* $.
\begin{itemize}
\item  The degree of an element
in $ S^k({\Gg}^*) \otimes \wedge^l {\Gg}^* $ is $ 2k+l$.
 With respect to this grading, $W({\Gg})$
is  clearly super-commutative.

\item  On generators, the differential  $d$ is defined by
\be
&& d (1\otimes \xi^i )\equal \xi^i \otimes
1-\half \sum_{j,k} f_{jk}^i \big(1\otimes (\xi^j \wedge
\xi^k)\big)\\
&& d (\xi^i \otimes 1)\equal  \sum_{j,k} f_{jk}^i (\xi^j \otimes \xi^k),
\ee
where $(\xi^i)_{i\equal 1}^{dim(\Gg)} $ is a basis of $\Gg^*$ and
$(f_{jk}^i)_{i,j,k\equal 1}^{dim(\Gg)} $ are the structure
constants of the Lie algebra $\Gg$ with respect to
the  dual basis of $(\xi^i)_{i\equal 1}^{dim(\Gg)}$.

\item The Lie group $G$ acts on $W(\Gg )$  by coadjoint action. 

\item For any $X\in \Gg$, the   contraction $ \ii_{X}$ is defined
to be an odd derivation, which, on the generators, is  given by
  $$\ii_{X}(\xi \otimes 1 )
\equal 0 \ \ \mbox{and } \   \ii_{X} (1 \otimes \xi )\equal \xi 
(X).$$
\end{itemize} 
In this case, the space of basic elements can be
identified with $S (\Gg^* )^{G}$.
\end{numex}

Another example is the following:

\begin{numex}
\label{ex:gspace}
Let $P$ be a $G$-space. Then the algebra $\Omega (P)$
 of differential forms on $P$ with the de Rham differential $d$,
the natural $G$-action and the contraction:
 $\ii_X \omega \equal \hat{X}\per \omega$,
where $\hat{X}$ is the infinitesimal vector field  on $P$
corresponding to $X\in \Gg$,
is a  $G$-differential algebra. In this case, the space of basic elements
can be identified with $\Omega (M)$.
\end{numex}

A {\em homomorphism} between two $G$-differential algebras $A$ and $B$
is an ${\mathbb N}$-graded differential algebra homomorphism
which commutes with the $G$-actions and the contractions. 
It is easy to see that the  $G$-differential algebras form
 a category.

\begin{defn}
\begin{enumerate}
\item A $G$-differential  simplicial algebra $ {\mathcal A}_\com$
 is a sequence of differential ${\mathbb N}$-graded
algebras $ (A_n)_{n \in {\mathbb N}} $
together with a sequence of morphisms of $G$-differential algebras, 
called \degeneracy  maps
$\epsilon_i^n: A_{n-1}~\to~A_{n}$, $i\equal 0,\dots,n$ and \face maps
$\eta_i^n:A_{n+1}~\to~A_{n} $, $i\equal 0,\dots,n$,  
which satisfy the co-simplicial relations:
$ \epsilon_j^{n}\epsilon_i^{n-1}
\equal   \epsilon_i^{n} \epsilon_{j-1}^{n-1}$ if $i<j$,
$ \eta_j^{n}\eta_i^{n+1} \equal  \eta_i^{n}\eta_{j+1}^{n+1} $ if $i \leq j$,
$ \eta_j^{n}\epsilon_i^{n+1} \equal  \epsilon_{i}^{n}  \eta_{j-1}^{n-1} $
if $i < j $, $ \eta_j^{n}\epsilon_i^{n+1} 
\equal   \epsilon_{i-1}^{n} \eta_{j}^{n-1}$
if $i > j+1 $ and $ \eta_j^{n}\epsilon_j^{n+1} 
\equal  \eta_j^{n}\epsilon_{j+1}^{n+1} \equal \id_{A_n} $.


\item A differential simplicial algebra is a $G$-differential simplicial
algebra where $G$ is the trivial group.
\end{enumerate}
\end{defn}

\begin{numex}
\label{ex:sim}
Given a ($G$-)simplicial manifold $M_\com\equal (M_n)_{n \in \nn}$, let
$A_n\equal  \Omega (M_n)$ be the differential algebra of
differential forms on $M_n$ equipped with the de Rham differential.
It is simple to see that $ (A_n)_{n \in {\mathbb N}} $
is a  ($G$-)differential  simplicial algebra with the
\degeneracy  and \face maps being the pull-back of
the \degeneracy  and \face maps of the simplicial manifold.  
\end{numex}

\begin{numex} \label{ex:4.6} In particular,
given a simplicial $G$-bundle  $Q_\com$ over $\gm_\com$,
according to Example \ref{ex:sim}, $(\Omega (Q_n ) )_{n \in {\mathbb N}}$,
is  a  $G$-differential  simplicial algebra.

Equivalently (see Proposition~\ref{prop:11corresp}),
if   $P$ is   a simplicial $G$-bundle over the groupoid
 $\Gamma\rightrightarrows \Gamma_0$, and  $Q\toto P$
is  the transformation groupoid,
then $Q_\com \equal ( Q_n )_{n \in \nn}$
is a principal $G$-bundle over the  simplicial  manifold
$\gm_\com$. Therefore one obtains  a
$G$-differential simplicial algebra $(\Omega (Q_n ) )_{n \in {\mathbb N}}$.
\end{numex}

The following construction gives us a useful  way
of constructing  $G$-differential simplicial algebras.

Let $A$ be a $G$-differential  unital algebra.
Set $A_n \equal A^{\otimes (n+1)}, \forall n\in \nn $ be the tensor  algebra
of graded algebras constructed according to  the Quillen rule \cite{Freed}.
The $G$-action on $A$ extends naturally to an 
action of $G$ on $A_n$ by the diagonal action. Similarly,
the contraction operation extends naturally to $A_n$ as
well  using the derivation rule. 
Define ${\eeta_i^{n}} : A^{\otimes n}
\to A^{\otimes n+1}$
and ${\eepsilon_i^n} : A^{\otimes n+2}\to A^{\otimes n+1}$ by
\be
&& {\eeta_i^n}(x_0 \otimes \dots \otimes x_{n-1})\equal 
x_0 \otimes \dots \otimes x_{i-1} \otimes 1 \otimes x_{i}
\otimes \dots \otimes x_{n-1} ,  \mbox{ and}\\
&& {\eepsilon_i^n}(x_0 \otimes \dots \otimes x_{n+1})\equal 
 x_0 \otimes \dots \otimes x_{i} x_{i+1}\otimes x_{i+2}
\otimes \dots \otimes x_{n+1}.  \  \ \\
\ee
It is simple to check that $\eeta_i^n$ and $\eepsilon_i^n$
satisfy all the compatibility conditions of \degeneracy maps
and \face maps. Indeed we have the following:

\begin{prop}
\label{pro:sim}
If $A$ is a $G$-differential  unital algebra, then
the sequence $A_\com: \equal (A^{\otimes (n+1)})_{n\in \nn} $
with the structures
described above is a $G$-differential simplicial algebra.
\end{prop}
\begin{pf}
The proof is a direct verification and is left for the reader.
\end{pf}

In particular, when $A$ is the Weil algebra $W(\Gg)$,
 one obtains the simplicial Weil algebra $W(\Gg)_\com$ (see \cite{Tondeur}):

\begin{defn}
 The simplicial Weil algebra $W({\Gg})_\com:
\equal  (W({\Gg})^{\otimes (n+1)})_{n\in \nn}$ is the $G$-differential
 simplicial algebra constructed as in Proposition \ref{pro:sim} when
$A$ is taken to be the Weil algebra $W({\Gg})$.
\end{defn}

One can also define  homomorphisms between (resp.   $G$-differential)
  simplicial algebras so that one obtains a category.

\begin{defn}
Let ${\mathcal A}_\com$ and ${\mathcal B}_\com$ be
differential simplicial algebras (resp. $G$-differential simplicial algebras).
 A homomorphism
is a sequence  of differential graded
algebra homomorphisms (resp. $G$-differential 
algebra homomorphisms) $ \phi_n:A_n \to B_n$ that  commute 
 with the \face and \degeneracy maps.
\end{defn}

\subsection{Cohomology and  fat realization of differential simplicial algebras}\label{sec:cup-prod}

Recall that the cohomology of an 
${\mathbb N}$-graded differential algebra $(A, d) $ is
$$H^* (A)\equal Ker\, d/ Im \, d,$$
and that $H(\cdot )$ is a covariant  functor from the
category of ${\mathbb N}$-graded differential algebras
to the category of graded algebras.
 
Similarly,  one can introduce 
 cohomology for simplicial differential algebras.
Let ${\mathcal A}_\com \equal  (A_n)_{n \in {\mathbb N}}$ be a
differential  simplicial algebra.
 By  $A_n^k $, we denote   the space of elements
of degree $k$ in $A_n$.  The collection  $ \{A_n^k|
(k,n) \in {\mathbb N}^2\}$ becomes
a  double complex with respect to the differential
$ d_n: A_n^k \to A_n^{k+1}$ of $A_n$ and the
map   $ \delta:A_n^k \to A_{n+1}^{k}$, which is defined to be
the alternate sum of the \degeneracy maps.
We denote by $H^* ( {\mathcal A}_\com ) $,  the cohomology
group  with respect to the
 total differential $ (-1)^n d +  \partial$.
Generalizing  Eq.~(\ref{eq:cup-productdR}),
there is a  cup-product
 $ \vee:\  A_{m}^{k} \otimes A_{n}^{l} \to A_{n+m}^{k+l}$  
 given by
\begin{equation}
 \label{eq:cupproduct} 
 a_1 \vee  a_2  \equal  (-1)^{kn}
\epsilon_{m+n}^{m+n} \smalcirc \cdots  \smalcirc \epsilon_{m+i}^{m+i} 
 \smalcirc \cdots   \smalcirc
\epsilon_{m+1}^{m+1}(a_1) \cdot \epsilon_{0}^{n+m}
 \smalcirc \cdots  \smalcirc \epsilon_{0}^{n+i} \smalcirc \cdots  \smalcirc
  \epsilon_{0}^{n+1}(a_2), 
\end{equation}
$\forall a_1 \in A_m^k, a_2 \in A_n^l $, that endows $H^*({\mathcal A}_\com) $
with a structure of super-commutative algebra. 

When ${\mathcal A}_\com $ is the simplicial
differential algebra of differential forms on
a simplicial manifold  $M_\com\equal (M_n)_{n \in \nn}$,
$H^*(A^*) $ reduces to the de Rham cohomology of the  simplicial manifold
as in Section \ref{section:deRham}.

\begin{prop}\label{prop:hfuncteur}
$H^*(\cdot )$ is a covariant  functor from the
category of simplicial  differential algebras
to the category of $\nn$-graded  algebras.
\end{prop}

There is an alternative way to
talk about the cohomology of a simplicial  differential algebra, namely 
via the  fat realization.
First, let us introduce the fat realization of differential
 simplicial algebras following the idea of Whitney, which
is  a standard construction for simplicial manifolds (see \cite{Dupont}).

Let $ \Delta_n$ be the $n$-simplex, {\em i.e}
$ \Delta_n\equal \{( t_0,\ldots,t_n)|\sum_{i\equal 0}^n t_i\equal 1,
t_0 \geq 0,\dots, t_n \geq 0 \}$. 
There are two sequences of maps $\tilde{\epsilon}_i^n :\Delta_{n-1} \to
\Delta_{n} $ and  $\tilde{\eta}_i^n :\Delta_{n+1} \to
\Delta_{n} $ defined by
\begin{eqnarray}
&&\tilde{\epsilon}_i^n (t_0,\dots ,t_{n-1})\equal 
(t_0,\dots,t_{i-1},0,t_{i},\dots,t_{n-1}), \label{eq:face1} \\
&&\tilde{\eta}_i^n (t_0,\dots ,t_{n+1})
\equal (t_0,\dots,t_i+t_{i+1},\dots,t_{n+1}).  \label{eq:deg1}
\end{eqnarray}

Assume that ${\mathcal A}_\com: \equal  (A_n)_{n \in {\mathbb N}}$ is a 
differential  simplicial algebra.  
Let
$$ \parallel {\mathcal A}_\com \parallel
\subset   \prod_{n \in {\mathbb N}}  \Omega ({\Delta_n}) \otimes A_n$$
be the subalgebra consisting of sequences
$ k_n \in   \Omega({\Delta_n} ) \otimes A_n $ such that
$\sup_{n\in {\mathbb N}} \mbox{deg}\,(k_n) < \infty$ and
\begin{eqnarray} 
&& ( \tilde{\epsilon}^{n*}_{i} \otimes \id) k_{n}\equal 
  ( \id \otimes \epsilon^{n}_i ) k_{n-1 }  , \ \ \ \ 
\mbox{and}  \label{eq:compatible} \\ 
 &&  ( \tilde{\eta}^{n-1*}_{i} \otimes \id) k_{n-1}\equal 
  ( \id \otimes \eta^{n-1}_i ) k_{n }  .\label{eq:compatible2}
\end{eqnarray} 
It is easy to see that
the differentials on $\Omega ({\Delta_n}) \otimes A_n , n\in \nn$,
respect the conditions (\ref{eq:compatible})-(\ref{eq:compatible2})
 so that $ \parallel {\mathcal A}_\com \parallel$
is a  differential algebra,
 called the {\em fat realization} of  $ {\mathcal A}_\com$.
The fat realization defines a covariant functor
from the category of differential simplicial
algebras to the category of differential algebras,
denoted by $\parallel \cdot \parallel$.
Similarly, the fat realization defines a covariant functor
from the category of $G$-differential simplicial
algebras to the category of $G$-differential algebras,
denoted by the same symbol $\parallel \cdot \parallel$.

\begin{numex}
When  ${\mathcal A}_\com $ is the simplicial
differential algebra of differential forms on
a simplicial manifold  $M_\com\equal (M_n)_{n \in \nn}$,
the fat realization $ \parallel{\mathcal A}_{\com}\parallel$
is (up to
some topological completion) the space of differential
forms on the geometric fat realization of $M_\com$ \cite{Dupont}.
More precisely, since $\Omega(\Delta_n) \otimes \Omega( M_n) $
is  a dense subset of $\Omega(\Delta_n \times  M_n) $,
 the fat realization $ \parallel {\mathcal A}_\com \parallel$
is a dense subset of the space of simplicial differential
forms on the geometric fat realization of $M_\com$.
\end{numex}

Given a simplicial differential algebra, we can apply the
fat realization functor to get a differential algebra
 and then the cohomology functor,
 or we can apply the cohomology functor directly.
 To compare these two approaches, we need the
following integration map following \cite{Dupont}:
\begin{eqnarray} \label{eq:above}
I: \parallel {\mathcal A}_\com  \parallel &\to & \bigoplus_n A_{n}\\
(k_n)_{n\in \nn}&\to & \sum_n  \int_{\Delta_n} k_n.
\end{eqnarray}

Note that the sum on the right-hand side of Eq. (\ref{eq:above})
is always  finite. Hence $I$ is well-defined. 

\begin{prop}\label{prop:proprieteintegration}
\begin{enumerate}
\item The integration map induces an isomorphism
$$I: H^* ( \parallel  {\mathcal A}_\com\parallel )\simeq  
H^* ( {\mathcal A}_\com ).$$ 
\item If $ \phi :{\mathcal A}_\com \to {\mathcal B}_\com$ is
a homomorphism of simplicial differential algebras,
then the following diagram commutes
$$
\comdia{ \parallel  {\mathcal A}_\com\parallel}{\parallel \phi \parallel}{ \parallel  {\mathcal B}_\com\parallel }{I}{}{I}{ {\mathcal A}_\com}{\phi}{ {\mathcal B}_\com}$$
\item Therefore one obtains a commutative diagram
 $$
\comdia{  H^* (\parallel  {\mathcal A}_\com\parallel)}{H(\parallel \phi \parallel)}{H^*( \parallel 
 {\mathcal B}_\com\parallel )}{I}{}{I}{ H^*({\mathcal A}_\com)}{H(\phi)}{ H^*({\mathcal B}_\com)}$$
\end{enumerate}
\end{prop}

\begin{pf} The proof is similar to that of Proposition 6.1
in \cite{Dupont} and we omit it here. 
\end{pf}

We summarize some important functorial properties of the above constructions
in the following

\begin{prop} \label{prop:fonctorialite}
\begin{enumerate}
\item The fat realization $\parallel \cdot \parallel$ is
a covariant functor from the category of differential
simplicial algebras (resp. the category of $G$-differential
simplicial algebras) to the category of differential
algebras (resp. the category of $G$-differential
 algebras). Both functors will be denoted by $R$.
\item Taking  basic elements    is
a covariant functor from the category of $G$-differential
simplicial algebras (resp. the category of $G$-differential
 algebras) to the category of differential simplicial
algebras (resp. the category of differential
 algebras). Both functors will be denoted by $B$.
\item The two functors $R\smalcirc B$ and
$B\smalcirc R$ are isomorphic as functors from
the category of $G$-differential
simplicial algebras  to  the category of differential
 algebras.
\end{enumerate}
\end{prop}

From now on, we will denote by
$ {\mathcal A}^{basic}_\com, \parallel {\mathcal A}_\com \parallel $
and $ H^*({\mathcal A}_{\com})$ the image of a $G$-differential simplicial
algebra ${\mathcal A}_{\com}$ under the functors $B$, $R$
and $H(\cdot)$ respectively.
 Accordingly, we will denote by $ \phi^{basic}, \parallel \phi \parallel $
and $H(\phi)$ the image of a homomorphism $\phi$ of $G$-differential
 simplicial algebras  under these functors.

\section{Chern-Weil map} \label{sec:CWM}

\subsection{Chern-Weil   map for $G$-differential simplicial algebras: first
construction} \label{subseq:Cher}

Following \cite{AM}, a {\em connection}  on a 
$G$-differential algebra $(A, d)$ is an
element $\theta \in (A^{1}\otimes \Gg)^G $ such that
for any $X\in \Gg$,
\begin{equation} \label{eq:defconnectionpart2}
\ii_X \theta \equal 1\otimes X.
\end{equation}
Note that  if $\phi : A \to B$ is a homomorphism of $G$-differential
algebras and  $\theta $ is a connection on $A$, then
$\phi  (\theta )\in (B^{1}\otimes \Gg)^G $  is
a connection on $B$.

 A connection induces a homomorphism of $G$-differential
algebras 
\begin{equation} \label{eq:defcwmapetape}c_{\theta}: W(\Gg) \to A ,
\end{equation}
which is defined on generators as follows. For any $ \xi \in \Gg^*$,
\be 
&&c_\theta (1\otimes \xi )\equal \la 1\otimes \xi , \theta \ra ,\\
&&c_\theta (\xi \otimes 1 )\equal \la 1\otimes \xi ,
d\theta +\half [\theta , \theta ]\ra.
\ee

Applying the basic functor, we get a morphism
 $z_{\theta}: S(\Gg^*)^G \to A^{basic}$,
where the differential on $S(\Gg^* )^G$
 is the zero differential. The map $z_\theta$ takes values in the 
subspace  $Z (A^{basic})$ of $A^{basic}$.
More precisely, we have
\begin{equation}\label{eq:explicitation} z_{\theta}(f)\equal f(\Omega),
   \ \ \ \ \forall f \in S(\Gg^*)^G \end{equation}
where $\Omega\equal d\theta +\half [\theta,\theta]\in (A^2 \otimes \Gg )^G$
 is the curvature of $\theta$.
 We call  $z_{\theta}$  the {\em  Chern-Weil map on the cochain level}
 \cite{guillemin}.

Applying the cohomology functor,
 one obtains a morphism $w_\theta : S(\Gg^*)^G \to H^*(A^{basic})$
called the  {\em  Chern-Weil map}.
The following result is well-known \cite{guillemin}.

\begin{prop}
 \label{prop:propbasealgweil}
\begin{enumerate}
\item The Chern-Weil map $w_{\theta}:S(\Gg^*)^G \to H^*(A^{basic}) $
 does not depend on the connection. 
In the sequel, we  denote $w_{\theta}$ by $ w$.  
\item If   $\phi : A \to B$ is a homomorphism of $G$-differential
algebras and  $\theta $ is a connection on $A$, then
$$c_{\phi (\theta )}\equal \phi \smalcirc c_\theta, \ \ 
z_{\phi (\theta )}\equal \phi \smalcirc z_\theta, \ \ \mbox{and }
w_{\phi (\theta )}\equal \phi \smalcirc w_\theta. $$
\end{enumerate}
\end{prop}

\begin{defn}
Let  ${\mathcal A}_\com:\equal  (A_n)_{n \in {\mathbb N}}$ be a
$G$-differential  simplicial algebra. We call a connection
on $A_0$  a pseudo-connection on ${\mathcal A}_\com$.
\end{defn}

Starting with  a pseudo-connection on ${\mathcal A}_\com$, we will
construct a connection on its fat realization
$\parallel {\mathcal A}_\com \parallel$.
Let us define $(n+1)$ maps from $A_0$ to $A_n$ by

   \be \label{eq:mapssimpliciaux} &&
p_0^n\equal \epsilon_n^n \smalcirc  \cdots \smalcirc \epsilon_i^{i}
\smalcirc\cdots \smalcirc \epsilon_1^1 \\
&& p_i^n\equal \epsilon_0^n \smalcirc \cdots
\smalcirc \epsilon_0^{n-i+1}  \smalcirc
\epsilon_{n-i}^{n-i} \smalcirc \cdots \smalcirc \epsilon_1^1,
\hspace{1cm}   1 \leq i < n \\
 && p_n^n\equal \epsilon_0^n \smalcirc  \cdots
\smalcirc\epsilon_0^{i}\smalcirc \cdots \smalcirc 
\epsilon_0^1 .
\ee

Now define, for any $\alpha \in A_0 $, 
  \begin{equation} \label{eq:defconnect} \alpha_n
\equal  \sum_{i\equal 0}^n t_i \otimes
 p_i^n (\alpha) \in \Omega^0(\Delta_n) \otimes  A_n .
 \end{equation}
The sequence $ \tilde{\alpha}:\equal  (\alpha_n )_{n \in \nn}$ satisfies 
the compatibility conditions
of Eqs.~(\ref{eq:compatible}-\ref{eq:compatible2}). Therefore
 $ \tilde{\alpha} \in \parallel {\mathcal A}_{\com}  \parallel   $.
In particular, if
$\theta \in  (A^{1}_0\otimes \Gg )^G $ is a pseudo-connection
on ${\mathcal A}_\com$, then  
\begin{equation}    \label{eq:defconnectexpl}  \tilde{\theta} \equal 
\big( \sum_{i\equal 0}^n t_i\otimes p_i^n (\theta) \big)_{n \in \nn}
\in  \parallel {\mathcal A}_{\com}  \parallel^1 \otimes \Gg .\end{equation}

\begin{prop} 
$\tilde{\theta} \in \parallel {\mathcal A}_\com  \parallel^1 \otimes \Gg$
defines a connection on the 
fat  realization $\parallel {\mathcal A}_\com \parallel$.
\end{prop}

\begin{pf} Since all the maps $p_i^n, i\equal 0,\dots,n$, are homomorphisms
of $G$-differential algebras,
$ p_i^n(\theta)\in (A_n^1 \otimes \Gg)^G, i\equal 0,\dots,n$,
is a connection on $ A_n$. Since $ \sum_{i\equal 0}^n t_i \equal 1 $,
hence $\theta_n\equal \sum_{i\equal 0}^n t_i p_i^n (\theta)$
is also a connection
on $  A_n $. 
In particular, for any $n \in \nn$, 
 $\theta_n $ is  $G$-invariant
 and therefore $ \tilde{\theta} $
is $G$-invariant as well.
  To show that  $\tilde{\theta}$ is a connection, it suffices  to check
Eq.~(\ref{eq:defconnectionpart2}) for any $ X \in \Gg$. Now
$$\begin{array}{ccc} \hat{X} \per  \tilde{\theta} &
\equal  & (  \hat{X} \per\theta_n )_{n \in \nn}  \\
                     &
\equal  & (  1_{\Omega(\Delta_n) \otimes  A_n} \otimes X )_{n \in \nn} \\
 &\equal  & (  1_{\Omega(\Delta_n) \otimes  A_n}  )_{n \in \nn} \otimes X \\
 &\equal  & 1_{\parallel {\mathcal A}_{\com}  \parallel } \otimes X .
 \end{array}$$
This completes the proof.
\end{pf}

According to Eq. (\ref{eq:defcwmapetape}),
 one obtains a   homomorphism of $G$-differential
algebras
$$ c_{\tilde{\theta}}:W(\Gg) \to \parallel {\mathcal A}_\com \parallel .$$
Applying the basic functor, one obtains a homomorphism  of
differential algebras
$$ c_{\tilde{\theta}}^{basic}:
S(\Gg^* )^G \to \parallel {\mathcal A}_\com \parallel^{basic}.$$
Using Proposition \ref{prop:fonctorialite}(3), we have
 $\parallel {\mathcal A}_\com \parallel^{basic} \simeq
\parallel {\mathcal A}_\com^{basic} \parallel$.
Therefore we obtain  a homomorphism,
 denoted again  by $c_{\tilde{\theta}}^{basic}$:  
\begin{equation} c_{\tilde{\theta}}^{basic}:
 S(\Gg^* )^G \to  \parallel {\mathcal A}_\com^{basic} \parallel .
\end{equation} 
Composing  $ c_{\tilde{\theta}}^{basic}$  with the cohomology functor
and using the isomorphism $H^* (S(\Gg^* )^G)\cong 
S(\Gg^* )^G$, one  obtains a
homomorphism of (graded) algebras 
$$ H(c_{\tilde{\theta}}^{basic}) :
S(\Gg^* )^G \to H^* (\parallel {\mathcal A}_\com^{basic} \parallel).$$
 Finally composing $H(c_{\tilde{\theta}}^{basic}) $  with
the integration map $I$, we obtain a map
\begin{equation}\label{eq:defCH-W}
\cw_\theta : S(\Gg^* )^G \to H^*  ( {\mathcal A}_\com^{basic}),
\end{equation}
which  is called the {\em Chern-Weil map}.

Note that $I: H^* (\parallel {\mathcal A}_\com^{basic} \parallel ) \to
H^*  ( {\mathcal A}_\com^{basic}) $ is a homomorphism of (graded) algebras.
Eq.~(\ref{eq:defCH-W})
implies that $\cw_\theta : S(\Gg^* )^G \to H^*  ( {\mathcal A}_\com^{basic}) $
is a homomorphism of algebras, 
whose image lies in $H^{even} ( {\mathcal A}_\com^{basic}) $.

On the cochain level,
composing  $ c_{\tilde{\theta}}^{basic}$  with the integration map
$I: \Omega (\parallel {\mathcal A}_\com^{basic} \parallel ) \to
\Omega  ( {\mathcal A}_\com^{basic}) $,
 we obtain a map $z_{\theta}$,
called the {\em  Chern-Weil map on the cochain level}:
   \begin{equation}\label{eq:defExplCH-Wcochain} 
z_{\theta} :  S(\Gg^* )^G \to Z  ( {\mathcal A}_\com^{basic}),
  \end{equation}
where $Z  ( {\mathcal A}_\com^{basic})$ denotes the space of cocycles 
in ${\mathcal A}_\com^{basic}$.

\begin{numex} \label{ex:simpl}
Consider the simplicial Weil algebra $W(\Gg )_\com $.
It is simple to check that
$\ttheta \equal \sum  (1\otimes \xi^i) \otimes X_i \in W(\Gg )^1 \otimes \Gg$
is a pseudo-connection, where $(\xi^i)_{i\equal 1}^{dim(\Gg)}$
and $(X_i)_{i\equal 1}^{dim(\Gg)}$  are dual basis
of $\Gg^*$ and $\Gg$. The  Chern-Weil map on the cochain level
induces a map 
\begin{equation}
\label{eq:Z}
Z:S(\Gg^* )^G \to Z^*(W(\Gg )_\com) 
\end{equation}        
and the Chern-Weil  map induces a 
morphism 
\begin{equation}
\label{eq:T}
T:  S(\Gg^* )^G  \to H^* (W(\Gg )_\com^{basic}).
\end{equation}
Indeed $T$ is an isomorphism according to Theorem 5.5   \cite{Tondeur}.
\end{numex}

\begin{them}
Let ${\mathcal A}_\com$ be a $G$-differential simplicial algebra.
Then
\begin{enumerate}
\item a pseudo-connection  $\theta \in  (A^{1}_0\otimes \Gg )^G $  
induces a canonical map
 \begin{equation}\label{eq:defExplCH-W}
z_{\theta} :  S(\Gg^* )^G \to Z  ( {\mathcal A}_\com^{basic}),
  \end{equation}
where $Z  ( {\mathcal A}_\com^{basic})$ denotes the space of cocycles
 in ${\mathcal A}_\com^{basic}$.
On the cohomology level, $z_{\theta}$ induces
a morphism  
\begin{equation}
\label{eq:defCH-W2}
\cw_\theta : S(\Gg^* )^G \to H^*  ( {\mathcal A}_\com^{basic}),
\end{equation}
which  is called the {\em Chern-Weil map}.                       
\item The Chern-Weil map $ w_{\theta}$ does not depend on the choice of
  pseudo-connections on ${\mathcal A}_\com$.
In the sequel, we will denote $w_{\theta}$ by $ w$. 
\end{enumerate}
\end{them}
\begin{pf} 1) is already proved.
2)  follows from Proposition \ref{prop:propbasealgweil}.
\end{pf}

\begin{numrk}
It would be interesting to compare our construction above 
with the  non commutative Chern-Weil map introduced
by Alekseev-Meinrenken \cite{AM03}.
\end{numrk}

\subsection{Second construction} \label{subseq:seco}

We now describe our second construction of the Chern-Weil
map.
Assume that ${\mathcal A}_\com:\equal  (A_n)_{n \in {\mathbb N}}$
is a $G$-differential  simplicial algebra and
$\theta \in  (A^{1}_0\otimes \Gg)^G $ is  a pseudo-connection
on ${\mathcal A}_\com$.

For any $0 \leq i \leq n$, $\theta_i^n \equal p_i^n (\theta)
\in (A_n^1 \otimes \Gg)^G $ is a connection on $A_n$. 
Thus we have 
$(n+1)$-homomorphisms of $G$-differential algebras
 $W(\Gg )\stackrel{c_{\theta_i^n}}{\lon} A_n $, $i\equal 0,\dots,n$.
Define $c: W(\Gg )_\com \to {\mathcal A}_\com$
by
\begin{equation}\label{eq:defphi}
c (x_0 \otimes \cdots \otimes x_n)\equal  c_{\theta_0^n} (x_0) \cdots
 c_{\theta_n^n} (x_n ), \ \ \ \forall x_i \in W(\Gg ), \ \ 
i\equal 0, \cdots , n.
\end{equation}

\begin{lem}
$c: W(\Gg )_\com \to {\mathcal A}_\com$ is a  homomorphism
 of $G$-differential simplicial algebras.
\end{lem}

\begin{pf}
The proof is a direct verification and is left to the reader. \end{pf}

Now applying the basic functor to Eq.~(\ref{eq:defphi}) yields
a homomorphism of differential simplicial algebras:
$$ c^{basic}: \ W(\Gg )_\com^{basic} \to {\mathcal A}_\com^{basic}.$$
Hence
\begin{equation}
\label{eq:c}
c^{basic}: Z^*(W(\Gg )_\com^{basic} )\to Z^*({\mathcal A}_\com^{basic}).
\end{equation}
By composing with the map $Z: S(\Gg^* )^G \to Z^* (W(\Gg )_\com^{basic} )$
 defined as in Eq. (\ref{eq:Z}), we obtain an homomorphism
\begin{equation}
\label{eq:zs}
z^s_{\theta} : S(\Gg^*)^G \to   Z^* ( {\mathcal A}_\com^{basic} ).
\end{equation}  
Applying the cohomology functor, we obtain a morphism of algebras
\begin{equation}
\label{eq:ws}
w^s_{\theta} : S(\Gg^* )^G \to  H^* ({\mathcal A}_\com^{basic} ). 
\end{equation}
We see in Corollary \ref{corr:construction2}  that indeed 
 $ z^s_{\theta} \equal  z_{\theta} $ and $ w^s_{\theta}\equal  w $.

\subsection{Properties of the Chern-Weil map} \label{subsec:prop}

\begin{prop}
\label{prop:prop} 
Let $\phi: {\mathcal A}_\com \to {\mathcal B}_\com$ be a homomorphism
of $G$-differential simplicial algebras,
 and $ \theta\in (A_0^1 \otimes \Gg)^G $ a pseudo-connection
on $ {\mathcal A}_\com$. Then the following diagrams commute
 \begin{equation}\label{eq:tri1}
\begin{array}{ccccccc}
S({\Gg}^*)^G  & \stackrel{z_{\theta}}{\to}&{Z^*}({\mathcal A}^{basic}_{\com}) & & S({\Gg}^*)^G &\stackrel{w}{\to}  & {H^*}({\mathcal A}_{\com}^{basic} )\\ 
&\sediag{z_{\phi(\theta)}} &\rdiag{\phi^{basic}} &\mbox{and}  &
 &\sediag{w} &\rdiag{H(\phi^{basic})} \\ 
& & {Z^*}({\mathcal B}^{basic}_{\com})& & & & {H^*}({\mathcal B}_{\com}^{basic} )\\ \end{array}
\end{equation}
\end{prop}
\begin{pf} From Eq.~(\ref{eq:defconnectexpl}), it follows that
$$ \begin{array}{ccc}
 \parallel \phi \parallel (\tilde{\theta})
& \equal  &\big( \phi( \sum_{i\equal 0}^n t_i \otimes
p_i^n (\theta) )\big)_{n \in \nn} \\ 
   &\equal
& \big(  \sum_{i\equal 0}^n t_i \otimes
p_i^n (\phi(\theta)) \big)_{n \in \nn}  \\ 
   &\equal  &    \tilde{\phi(\theta)}.  \hspace{3cm} \\
\end{array} $$

According to Proposition \ref{prop:propbasealgweil},
we have the commutative diagram
\begin{equation} \label{eq:ajout1}
\comtri{ W({\Gg})}{c_{\tilde{\theta}}} {  \parallel  {\mathcal A}_{\com}  \parallel}{c_{\tilde{\phi(\theta)}}}{\parallel \phi \parallel}{ \parallel  {\mathcal B}_{\com}  \parallel}\end{equation}

Applying the basic functor $B$, one  obtains the commutative
diagram:
$$\comtri{S({\Gg}^*)^G }{c_{\tilde{\theta}}^{basic}}{ \parallel  {\mathcal A}_{\com}  \parallel^{basic}}{ c^{basic}_{ \tilde{\phi(\theta)} }}{ \parallel \phi \parallel^{basic}}{  \parallel  {\mathcal B}_{\com}   \parallel^{basic}  } . $$

By Proposition \ref{prop:fonctorialite}, the realization and basic functors commute.  Therefore, we have
\begin{equation}\label{eq:relationbasic}\comtri{S({\Gg}^*)^G }{c_{\tilde{\theta}}^{basic}}{ \parallel  {\mathcal A}_{\com}^{basic}  \parallel}{ c^{basic}_{ \tilde{\phi(\theta)} }}{ \parallel \phi^{basic} \parallel}{  \parallel  {\mathcal B}_{\com}^{basic}   \parallel  }  \end{equation}
According to Proposition \ref{prop:proprieteintegration},
 $  I \smalcirc \parallel \phi^{basic} \parallel
 \equal  \phi^{basic} \smalcirc I $.
 Composing    with the integration map, we have 

$$\begin{array}{ccccc}
 S(\Gg^*)^G &\stackrel{c_{\tilde{\theta}}^{basic}}{\to} &
\parallel  {\mathcal A}_{\com}^{basic}  \parallel & \stackrel{I}{\to} &{\mathcal A}_{\com}^{basic} \\
  & \sediag{  c_{\tilde{\phi(\theta)}}^{basic}   } &\rdiag{ \parallel \phi^{basic} \parallel }& &\rdiag{\phi^{basic}}  \\
  &         & \parallel  {\mathcal B}_{\com}^{basic}  \parallel&\stackrel{I}{\to} & {\mathcal B}_{\com}^{basic} \\
 \end{array}  $$
Since $I \smalcirc c_{\tilde{\theta}}^{basic}$ and $ I \smalcirc c_{\tilde{\phi(\theta)}}^{basic} $ 
take values in $Z( {\mathcal A}_{\com}^{basic})$
and $Z( {\mathcal B}_{\com}^{basic}) $ respectively,
one obtains the commutative diagrams  (\ref{eq:tri1}).
\end{pf}

\begin{cor}\label{corr:construction2} 
For any $G$-differential simplicial algebra $ {\mathcal A}_\com $ and any 
pseudo-connection $\theta$ on $ {\mathcal A}_\com $, we have
$z_{\theta} \equal z^s_{\theta}$ and $ w \equal w^s$.
Thus  the first  and the second Chern-Weil  constructions 
coincide.
\end{cor}
\begin{pf} 
Let $\eta \in W(\Gg)^1\otimes \Gg$
be   the pseudo-connection constructed as in Example \ref{ex:simpl}.
It is simple to  check that
$ c(\eta)\equal \theta $,
 where $c: \ W(\Gg)_\com \to {\mathcal A}_\com$
is the map  defined by Eq.~(\ref{eq:defphi}).
By Proposition \ref{prop:prop}, we have $c^{basic} \smalcirc z^s_{\theta}\equal z_\theta$.
Passing to the cohomology, we have $w_\theta^s\equal w$.
  \end{pf}

\begin{rmk}
From Eq.~(\ref{eq:ajout1}),  it follows that
for any  homomorphism of $G$-differential simplicial algebras 
$\phi: {\mathcal A}'_\com  \to {\mathcal A}_\com$, the equality
$c_{\tilde{\phi({\theta})}} \equal \parallel \phi \parallel
\smalcirc c_{\tilde{\theta}}
$ holds.
Composing with $I$ we obtain
  $$  I \smalcirc c_{\tilde{\phi(\theta)}} \equal  I \smalcirc \parallel \phi
\parallel \smalcirc c_{\tilde{\theta}} .$$
Therefore, by Proposition \ref{prop:proprieteintegration}(2), we obtain
  \begin{equation}\label{eq:ajout2}
I \smalcirc c_{\tilde{\phi(\theta)}}  \equal 
\phi \smalcirc  I \smalcirc  c_{\tilde{\theta}}. \end{equation}
In particular, if ${\mathcal A}'_\com \equal W(\Gg)_\com$
is endowed with the pseudo-connection $\eta \in W(\Gg)^1\otimes \Gg$
as in Example \ref{ex:simpl}, $\zeta$
is a pseudo-connection on ${\mathcal A}_\com $
and $c: W(\Gg)_\com\to {\mathcal A}_\com $ 
is the homomorphism of $G$-differential simplicial algebras
constructed as in
Eq.~(\ref{eq:defphi}), then $c(\eta )=\zeta$. Thus Eq.~(\ref{eq:ajout2})
implies  that
  \begin{equation}\label{eq:withcs}  I \smalcirc c_{\tilde{\zeta}}
\equal 
c \smalcirc  I \smalcirc  c_{\tilde{\eta}}.
\end{equation}
\end{rmk}

\section{Chern-Weil map for principal $G$-bundles over groupoids}
\label{sec:CWgroupoides}

\subsection{Main theorem}

In this subsection, we apply the results of
Section \ref{sec:CWM} to the case of
a principal $G$-bundle over a groupoid. Let $ P \stackrel{\pi}{\to} \gm_0$
be a principal $G$-bundle over $\gm \toto \gm_0$.
Then, according to Example \ref{ex:4.6}, $ \Omega(Q_{\com})$
 is a $G$-differential 
simplicial algebra and a pseudo-connection $\theta \in \Omega^1(P) \otimes \Gg$
defines a pseudo-connection on the $G$-differential
simplicial algebra $\Omega(Q_{\com})$,
where $Q \toto P$ is the transformation groupoid.
Note that $\Omega(Q_{\com})^{basic}\equal \Omega(\gm_{\com}) $.
Therefore, one obtains a Chern-Weil map
  $w_P: S(\Gg^*)^G \to H^*_{dR}(\gm_\com) $ and a Chern-Weil map
$z_\theta: S(\Gg^*)^G \to Z^*_{dR}(\gm_\com) $ 
on the cochain level.

\begin{them}\label{th:grou6.1}
\begin{enumerate}
\item  Associated to any pseudo-connection
$\theta \in \Omega^1 (P)\otimes \Gg$, there is a canonical
map
$$z_\theta: S(\Gg^*)^G \to Z^*_{dR}(\gm_\com), $$
called the  Chern-Weil map on  the cochain level, where
$Z^*_{dR}(\gm_\com)$ is the space of closed forms.
On the cohomology level, $z_\theta$ induces an algebra
homomorphism
$$w_\theta: S(\Gg^*)^G \to H^*_{dR}(\Gamma_\com), $$
which is independent of  the choice of  pseudo-connections and is
denoted by $w_P$.
Moreover,  $z_\theta$ is completely determined
by the total pseudo-curvature (Proposition \ref{pro:5.11}).
\item If $ \phi$ is a strict homomorphism
from $ \gm' \toto \gm_0'$ to $ \gm \toto \gm_0$,
then the following diagrams $$ \begin{array}{ccccccc}
 S({\Gg}^*)^G &\stackrel{z_{\phi^* \theta}}{\to} & Z_{dR}^*(\gm_{\com}') &
&  S
({\Gg}^*)^G & \stackrel{w_{P'}}{\to} & H^*_{dR}(\gm_{\com}')   \\
  & \stackrel{\searrow}{z_{ \theta}}&
\uparrow \phi^* &\mbox{and} & & \stackrel{
\searrow}{w_P}&  \uparrow \phi^* \\
  & &  Z_{dR}^*(\gm_{\com}) & & & & H^*_{dR}(\gm_{\com})
   \end{array}$$
commute, where $P'\equal  \phi^* P$ is the pull-back of $P$ via $\phi$.
\end{enumerate}
\end{them}                                  

The following proposition  lists some important properties
of this Chern-Weil map. 

\begin{them}\label{th:grou1}
\begin{enumerate}
\item If $ \Gamma$ is a manifold $M$, then $w_P: S(\Gg^*)^G \to H^*_{dR}(M)
$
reduces to the usual Chern-Weil map \cite{guillemin}.
 \item If $ \Gamma$ is a Lie group $G$,  $P$ is
 the $G$-bundle  $G \to \cdot $ and $\theta$ is the left
Maurer-Cartan form, then  $z_\theta: S(\Gg^*)^G \to Z_{dR}^*(G_\com) $
and $w_P: S(\Gg^*)^G \to H^*_{dR}(G_\com) $
coincide with the Bott-Shulman maps \cite{Bott,BSS}.
\end{enumerate}
\end{them}               
\begin{pf} 1) When $\Gamma \equal M$, the one-forms $\alpha_n$ defined by 
Eq.~(\ref{eq:mapssimpliciaux}) are all equal. Therefore, 
 $  \tilde{\theta} \equal  (1_{\Omega(\Delta_n)} \otimes \theta)_{n \in \nn}$
and 
 \begin{equation}\label{eq:suiteannule} c_{\tilde{\theta}}^{basic}
\equal    (1_{\Omega(\Delta_n)}  )_{n \in \nn}
\otimes c_{\theta}^{basic} . \end{equation}
The integration map  $I$ is equal to zero on $\Omega^0(\Delta_n)
\otimes \Omega^*(M)$ unless $n \equal 0$. Composing
Eq~.(\ref{eq:suiteannule}) with $I$, we obtain
 $  z_{\theta} \equal  I \smalcirc c_{\tilde{\theta}}^{basic}
\equal  c_{\theta}^{basic}$.
The conclusion follows by passing to the cohomology.

2) Let $\theta \in \Omega^1 (G)\otimes \Gg$
 be the left Maurer-Cartan form of $G$. Then our construction 
 reduces to  that as in \cite[Chapter 6]{Dupont}.
 Hence, the conclusion  follows from  \cite{BSS}.
\end{pf}

Now we turn to the study of  the relation between
the Chern-Weil map and the total pseudo-curvature.

As above,  let $\theta\in \Omega^1 (P)\otimes \Gg$ be 
a pseudo-connection and $\Omega_{total}\equal \partial \theta+
\Omega\in \Omega^1 (Q)\otimes \Gg \oplus \Omega^2 (P)\otimes \Gg$
its pseudo-curvature. Let  $ {\mathcal C} $ be the subalgebra
of $ \Omega(Q_{\com})$  generated (under the cup-product) by the
images of  both maps $\wedge \Gg^* \to \Omega^1 (Q)$
and $S(\Gg^*) \to \Omega^2 (P)$ induced by $\partial \theta$ and
$\Omega$ respectively.  Let $ {\mathcal D} \subset  \Omega(\gm_{\com})$
be the subalgebra of $\Omega(\gm_{\com})$ consisting of
basic elements of $ {\mathcal C} $.

%

\begin{prop}
\label{pro:5.11} 
For  any $ f \in S(\Gg^* )^G$, $z_\theta(f)$ is a
non-commutative polynomial $P_f(\del\theta,\Omega)$ in $\del\theta$
and $\Omega$. In particular, $z_{\theta}$ belongs to ${\mathcal D}$.
\end{prop}

In other words the Chern-Weil map is completely determined 
by the total pseudo-curvature.
\begin{pf}
The curvature of $ \tilde{\theta}$ on $ \parallel Q_\com \parallel $
is the simplicial $\Gg$-valued
$2$-form  $ \tilde{\Omega}\equal d\tilde{\theta}+ \half [\tilde{\theta},
\tilde{\theta}]$. More precisely
\begin{equation} \label{eq:explicitecourbure}  \tilde{\Omega}_n \equal 
 \sum_{i\equal 0}^n d t_i  \wedge  \theta_i +
\sum_{i\equal 0}^n  t_i  \Omega_i -
 \frac{1}{2}  \sum_{i\equal 0}^n t_i [\theta_i,\theta_i]  +
 \frac{1}{2}  \sum_{i\equal 0}^n \sum_{j\equal 0}^n t_i  t_j
[\theta_i , \theta_j].
 \end{equation}
where   $  \theta_i\equal (p^n_i)^* \theta \in \Omega^1 (Q_n )\otimes
\Gg$ and  $\Omega_i\equal (p^n_i)^* \Omega \in \Omega^2 (Q_n ) \otimes \Gg $.

Introducing new variables  $\eta_i\equal \theta_{i+1}-\theta_i$, $i=1,\cdots,n$
we have
\begin{eqnarray*}
\lefteqn{\sum_{i\equal 0}^n t_i [\theta_i,\theta_i]  -
 \sum_{i\equal 0}^n \sum_{j\equal 0}^n t_i  t_j[\theta_i,\theta_i]}\\
&\equal &\sum_{i\equal 0}^nt_i
[\theta_0+\sum_{k<i}\eta_k,\theta_0+\sum_{l<i}\eta_l]
-\sum_{i,j}t_it_j[\theta_0+\sum_{k<i}\eta_k,\theta_0+\sum_{l<j}\eta_l]\\
&\equal &\sum_{k\equal 0}^{n-1}
\sum_{l\equal 0}^{n-1} \big( \sum_{i>max(k,l)}t_i \big) [\eta_k,\eta_l]
-\sum_{k\equal 0}^{n-1}   \sum_{l\equal 0}^{n-1}\sum_{i>k,j>l}
t_it_j  [\eta_k,\eta_l]
\end{eqnarray*}

Moreover the identity $\sum_{i \equal 0}^n dt_i \wedge \theta_i
\equal  -\sum_{i=0}^{n-1} ds_i \wedge \phi_i $ holds
where $s_i  \equal  \sum_{k \equal 0}^i t_k$ for $i=0,\cdots,n-1$. 
Therefore
Eq.~(\ref{eq:explicitecourbure})  can be rewritten as

\begin{equation} \label{eq:real72} \begin{array}{rcl}     \tilde{\Omega}_n&
\equal & -\sum_{i\equal 0}^{n-1} d
 s_i \land (\theta_{i+1}-\theta_{i}) +
\sum_{i\equal 0}^n t_i \Omega_i \\ & &
-\frac{1}{2}
\sum_{k\equal 0}^{n-1}\sum_{l\equal 0}^{n-1}
\big(\big( \sum_{i>max(k,l)}t_i\big) 
+ \big( \sum_{i>k,j>l}
t_it_j \big)\big)  [\theta_{k+1}-\theta_k,\theta_{l+1}-\theta_l]. \\
\end{array} \end{equation}

For any $f \in S(\Gg^* )^G$, we have $ p^*(z_\theta (f))\equal I 
\smalcirc f(\tilde{\Omega}) $,
where $p^*: \Omega(\gm_\com)\to \Omega(Q_\com)$ is the pull-back
of the projection   $Q \to \gm $.
Now, by Eq.~(\ref{eq:real72}),
it is simple to see that
 $f(\tilde{\Omega}_n) $ is a linear combination
of terms of the form:
   $$ \omega \wedge  f_0 (\Omega_0) \wedge g_0 (\theta_1 -\theta_0 )
\wedge \cdots \wedge  f_i (\Omega_i) \wedge g_i (\theta_{i+1} -\theta_i )
\wedge\cdots \wedge f_{n-1} (\Omega_{n-1}) \wedge g_{n-1} 
(\theta_{n} -\theta_{n-1} )\wedge f_n (\Omega_n ),$$
where $ \omega \in \Omega (\Delta_n) $,
$f_i\in S(\Gg^*), \ i\equal 0, \cdots , n$ and 
$g_i\in \wedge \Gg^*, \ i\equal 0, \cdots , n-1$.

Note that by applying the integration map to the above
forms, $\omega$ is integrated and    the   remaining 
terms  are linear combinations of the form
$ f_0 (\Omega_0) \wedge g_0 (\theta_1 -\theta_0 )
\wedge \cdots \wedge  f_i (\Omega_i) \wedge g_i (\theta_{i+1} -\theta_i )
\wedge\cdots \wedge f_{n-1} (\Omega_{n-1}) \wedge g_{n-1} (\theta_{n} -
\theta_{n-1} )\wedge f_n (\Omega_n )$. 
Now by the definition of cup-product, one can 
rewrite the latter as
\be && f_0 (\Omega_0) \wedge g_0 (\theta_1 -\theta_0 )
\wedge \cdots \wedge  f_i (\Omega_i) \wedge g_i (\theta_{i+1} -\theta_i )
\wedge\cdots \wedge f_{n-1} (\Omega_{n-1}) \wedge g_{n-1}
 (\theta_{n} -\theta_{n-1} )\wedge f_n (\Omega_n )
\\&& \equal f_0 (\Omega )\vee  g_0  (\partial \theta ) 
\vee \cdots \vee  f_i (\Omega) \vee g_i (\partial \theta )
\vee \cdots \vee  f_{n-1} (\Omega) \vee  g_{n-1} (\partial \theta)
\vee  f_n (\Omega). \ee
This completes the proof.
\end{pf}

 Consider the subcomplex $(\Omega^{\bullet} (\gm_0)^\gm, d)$ of the
de Rham complex $(\Omega^{\bullet} (\gm_0), d)$, where
  $\Omega^{\bullet} (\gm_0)^\gm \equal  \{\omega \in \Omega^{\bullet} (\gm_0)| \ 
\partial \omega \equal 0 \}$.
Let $ H^*(\gm_0)^\gm $ be its   cohomology group.
The natural inclusion $ i:\Omega^{\bullet} (\gm_0)^\gm \to \Omega^{\bullet} (\gm_{\com}) $
is a chain map and induces a homomorphism $ H^*(\gm_0)^\gm  
\to H^*_{dR}(\gm_\com)$.

When a connection exists, the Chern-Weil map  admits an explicit
simple form.

 \begin{them}\label{th:connectionchernweil}
Assume that the principal $G$-bundle  $ P \stackrel{\pi}{\to} \gm_0$ over
 the groupoid $\gm \toto \gm_0$ admits a connection $\theta
\in \Omega^1 (P) \otimes \Gg$. Then
  the  following diagrams  commute
$$\comtri{ S(\Gg^*)^G }{ z }{Z^*(\gm_0)^\gm }{z_\theta}{i}{Z^*(\gm_{\com})}
$$
and
 $$\comtri{ S(\Gg^*)^G}{ w}{ H^*_{dR}(\gm_0 )^{\gm }}{w_P}{  i }{H^*_{dR}(\gm_{\com})} $$
where  $z: S(\Gg^*)^G \to Z^*(\gm_0)^{\gm}$, and
$w: S(\Gg^*)^G \to H^*_{dR}(\gm_0)^{\gm}$
denote  the usual Chern-Weil maps
by forgetting about the groupoid action {\em i.e.},
 $z(f) \equal f(\Omega)\in Z^* ( P )^{basic}\cong
Z^* (\gm_0) $ and $w(f)=[z(f)] \in H^*(\Gamma_0)$, 
$\forall f \in S(\Gg^*)^G$, where
$\Omega \in \Omega^2 (P)\otimes \Gg$ is the curvature form.
\end{them}
\begin{pf}
 We have 
$$  \partial \Omega\equal  \partial (d\theta + \half[\theta,\theta])
\equal d \partial\theta + \half \big( [\partial\theta,\theta]+
[\theta,\partial \theta]\big) \equal 0 . $$
Hence the image of   $z: S(\Gg^*)^G \to Z^*(\gm_0)$  lies in
$Z^*(\gm_0)^\gm$. Therefore, the image of   $w : S(\Gg^*)^G \to H^*_{dR}(\gm_0)$
lies in  $ H^*_{dR}(\gm_0)^\gm$.

Since $\partial \theta \equal 0$, it follows that all
the terms $\theta_i -\theta_{i-1}$
in  Eq.~(\ref{eq:real72})  vanish for  $i \in \{0,\dots,n-1\}$.
Thus we have $\tilde{\Omega_n} \equal \sum_{i\equal 0}^n t_i \Omega_i$.
This  implies that for any $f \in S(\Gg^*)^G $,
 $I ( f(\tilde{\Omega_n}))\equal 0$ except for $n\equal 0$.
Therefore
the identities $z_\theta(f)
\equal I\big(f(\tilde{\Omega}_0)\big)\equal i(f(\Omega))$ hold.
\end{pf}

\begin{rmk}
 In the case of an \'etale  groupoid or the holonomy groupoid
of a foliation, the Chern-Weil map has been  studied by
Crainic and Moerdijk in \cite{Moerdjick}.
\end{rmk}

\begin{cor}
If a principal $G$-bundle  $ P \stackrel{\pi}{\to} \gm_0$ over
 a  groupoid $\gm \toto \gm_0$ admits a flat connection,
then the Chern-Weil map vanishes (except for degree $0$).
\end{cor}

\begin{numex}
The Chern-Weil map associated to a $G$-bundle 
over a finite group is always zero in strictly positive degree.
This is due to the fact that
this bundle admits a flat connection.
\end{numex}

\subsection{Universal Chern-Weil map}

According to Proposition \ref{prop:11corr}, a principal $G$-bundle
$P \to \gm_0$ over $\gm \toto \gm_0$ 
corresponds to
 a generalized homomorphism $\phi_P$
from $ \gm \toto \gm_0$ to $G \toto \cdot$.
  According to Eq.  (\ref{eq:defpullback}),
 such a generalized homomorphism induces a map
$ \phi^*: H^*_{dR}(G_{\com}) \to H^*_{dR}(\gm_{\com}) $.
Let  $ S(\Gg^*)^G \to  H^*_{dR}(G_{\com}) $ be the Bott-Shulman
map.
Composing these two morphisms, we obtain a morphism
    $$ w_P^{u} :  S(\Gg^*)^G \to H^*_{dR}(\gm_{\com}) $$
called the {\em universal Chern-Weil map}.
It is immediate from the definition that if $\phi: \gm' \to \gm$
is a generalized homomorphism then
\begin{equation}\label{eq:comportuniv} w^u_{P'}\equal
\phi^* \smalcirc w_P^u  \end{equation}
where $P'\equal \phi^* P$ is the pull-back of
$P \stackrel{\pi}{\to} \gm_0$ by $\phi$
as constructed in Proposition \ref{prop:pullback} and
$\phi^*:H^*_{dR}(\gm_\com') \to H^*_{dR}(\gm_\com)$
 is the map constructed as in Eq.~(\ref{eq:defpullback}).  

\begin{them}
\label{prop:coher} The universal Chern-Weil map and the
Chern-Weil map are equal. \end{them}
\begin{pf} According to
Theorem \ref{th:grou1}(2), we know that
the claim holds for
 the principal $G $-bundle $G \to \cdot$ over $G \toto \cdot $.

 Let $ Q'$ be the groupoid $Q \times G \toto P$ as in
 Eq.~(\ref{eq:productQ'}). And let
$pr:Q'\to G$ and $p:Q'\to \Gamma$ be the strict homomorphisms
 given by the natural projections on $G$ and $\gm$, respectively.
It is simple to see that the pull-back of the principal $G$-bundle
 $G \to \cdot$ by  $pr$ is isomorphic to
$ \tilde{P}\equal P \times G  \to P$. 
According to Eq.   (\ref{eq:comportuniv}) and
Theorem \ref{th:grou6.1}(2), we have
  \begin{equation} \label{eq:premiereetape}
   w^u_{\tilde{P}}\equal   w_{\tilde{P}} : S(\Gg^* )^G\to H^*_{dR}(Q'_\com ) .
\end{equation}
On the other hand, the pull-back of $P \to \gm_0 $ via
$p$ is also isomorphic to $ \tilde{P}\to P$. Hence
Eq.  (\ref{eq:comportuniv}) implies that 
$  w_{\tilde{P}}^u  \equal p^*  w_{P}^u $, while 
Theorem \ref{th:grou6.1}(2) implies
that $  w_{\tilde{P}}\equal  p^*  \smalcirc w_{\tilde{P}}. $ Therefore
   \begin{equation}\label{eq:avantiso}
 p^*  w^u_{P}\equal  p^*   w_P.
 \end{equation}

It is easy to see that $ Q'  \toto P$ is
indeed the pull-back of  $\gm \toto \gm_0$
via $P \stackrel{\pi}{\to} \gm_0 $. So these
two groupoids  are Morita equivalent,
 and therefore, according to Lemma \ref{lem:moritamap},
 $p^*: H^*_{dR}(\gm_{\com}) \to H^*_{dR}(Q_{\com}')$
is an isomorphism. Therefore,
 Eq.~(\ref{eq:avantiso}) implies that $  w^u_{P}\equal   w_P .$
\end{pf}

The following corollary follows immediately
from  Proposition \ref{prop:coher} and Eq.  (\ref{eq:comportuniv}).

\begin{cor}\label{cor:univgenehomo} Let $\phi$ be a generalized homomorphism from $\gm' \toto \gm_0'$
to $\gm \toto \gm_0$, then the following diagram commutes

$$\comtriup{S(\Gg^*)^G}{w_{\phi^* P}}{H^*_{dR}(\gm'_{\com})}{w_{P}}{\phi^*}{ H^*_{dR}(\gm_{\com} )} $$
where $P'\equal \phi^* P $ is the pull-back of $P$ via $\phi$ and $\phi^*: H^*_{dR}(\gm_{\com}) \to  H^*_{dR}(\gm_{\com}')$ is the morphism 
induced  by $\phi$ as given by Eq.~(\ref{eq:defpullback}).
\end{cor}
\begin{rmk}
From Corollary \ref{cor:univgenehomo}, it follows
that the universal Chern-Weil map is defined
for $G$-bundles over differential stacks.
Note that generalized homomorphisms correspond to smooth maps between their
stacks. Theorem \ref{prop:coher} means that, when
 a presentation of a differential stack, {\em i.e.} a Lie groupoid, is chosen,
the universal Chern-Weil map can be computed using a pseudo-connection,
as in the manifold case.
\end{rmk}

\subsection{Equivariant principal $G$-bundles}

For a compact Lie group $H$ and a smooth manifold $M$ on which $H$ acts,
 there are many different
ways to define the equivariant cohomology $H^*_H(M)$.
One definition, called the simplicial model, is to define
$H^*_H(M)$ as  $ H^*_{dR}(\Gamma_\com)$,
 where  $\gm$
is the transformation groupoid  $ H \times M \toto M$. 
Another  definition, called the Weil model, is to define $H_H^*(M) $
as the cohomology of the complex
$ \big(W(\Hh) \otimes \Omega (M)\big)^{H-basic}$, 
where $\Hh$ is the Lie algebra of $H$.
These two models are known to be equivalent (see Section 4.4 \cite{guillemin}).

As an application of the tools developed in the previous sections,
 in  what follows,
we describe an explicit  chain map that establishes
a  quasi-isomorphism from the Weil model
to the simplicial model.
Consider the right principal $H$-bundle $H \to \cdot$
 over the groupoid $H \toto \cdot $ together
with the pseudo-connection $\zeta$ equal to the Maurer-Cartan form,
i.e., 
 $\zeta\equal \sum_{i\equal 1}^{dim(\Hh)} \zeta^i \otimes X_i$ where $\{X_1,
\dots,X_{dim(\Hh)}\}$ is a basis of $\Hh$ 
and  $\zeta^1,\dots,\zeta^{dim(\Hh)}\in \Omega^1(H)$ be the dual basis
of the left invariant vector fields associated to $X_1,
\dots,X_{dim(\Hh)}\}$.  As in Example \ref{ex:simpl},
 let  $\eta\equal  \sum_i ( 1 \otimes \xi^i) \otimes X_i$
be a connection on the algebra $W(\Hh) $,
where $\{\xi^1,\dots,\xi^{dim(\Hh)} \}\in \Hh^*$
is the dual basis of    $\{ X_1,\dots,X_{dim(\Hh)} \}$.

Note that we can  consider $\eta$ as a pseudo-connection on $ W(\Hh)_\com$.
Similarly, let $H \times H \toto H$ be the pair groupoid.
We can  consider $\zeta$ as a pseudo-connection on $\Omega\big((H \times H)_\com\big) $.
Let $c:W(\Hh)_\com \to \Omega\big((H \times H)_\com \big)$ be the homomorphism
of $H$-differential simplicial algebras
 defined as in Eq.~(\ref{eq:defphi})
using the pseudo-connection $\zeta$.
It is clear that $c(\eta) \equal \zeta $. Hence,
Eq.~(\ref{eq:withcs}), applied to the case where ${\mathcal A}_\com \equal  \Omega\big((H \times H)_\com $
and $\theta \equal \eta $
implies that the following two  sequences of compositions
of chain maps are equal
  $$ W(\Hh) \stackrel{c_{\tilde{\zeta}}}{\to} \parallel \Omega\big((H \times H)_\com\big) \parallel \stackrel{I}{\to} \Omega\big((H \times H)_\com\big) $$
and 
  $$  W(\Hh) \stackrel{c_{\tilde{\eta}}}{\to} \parallel W(\Hh)_{\com} \parallel \stackrel{I}{\to} W(\Hh)_{\com} \stackrel{c}{\to} \Omega
((H \times H)_\com).   $$

 In other words, we have a commutative diagram
$$
\xymatrix{ \parallel \Omega\big((H \times H)_\com\big) \parallel      
\ar[r]^I&   \Omega\big((H \times H)_\com\big)      \\
 W(\Hh) \ar[u]\ar[d]_{c_{\tilde{\eta}}} &  \\
 \parallel W(\Hh)_{\com} \parallel  \ar[r]^I &  W(\Hh)_{\com}    \ar[uu]^c  \\
}
$$
Taking the tensor product with $ \Omega(M)$, we obtain 
\begin{equation}\label{eq:tableau}
\xymatrix{
 \parallel \Omega\big((H \times H)_\com\big) \parallel  \otimes \Omega(M)     \ar[r]^I& 
 \parallel \Omega\big((H \times H)_\com\big) \parallel  \otimes \Omega(M)   \otimes \Omega(M)    \\
W(\Hh)  \otimes \Omega(M)  \ar[u]   \ar[d]_F &  \\ 
 \parallel W(\Hh)_{\com} \parallel  \otimes \Omega(M)  \ar[r]^I &  
   W(\Hh)_{\com} \otimes \Omega(M) \ar[uu]^{ c \otimes \id }  
}
\end{equation}
where $F\equal  c_{\tilde{\eta}}\otimes \id$.
Let $ R \toto H \times M $ be the product of the pair 
groupoid $H \times H \toto H$ with the manifold
$ M \toto M$. 
Since we have the inclusion $\Omega((H \times H)_\com)\otimes \Omega(M)
\subset \Omega(R_\com )   $ and $   \parallel \Omega\big((H \times
H)_\com\big)\parallel\otimes \Omega(M)  \subset \parallel \Omega(R_\com )
\parallel$, we have 
\begin{equation}
\xymatrix{
\parallel \Omega(R_\com) \parallel  \ar[r]^I&  \Omega(R_\com)    \\
W(\Hh)  \otimes \Omega(M)  \ar[u]   \ar[d]_F &  \\ 
 \parallel W(\Hh)_{\com} \parallel  \otimes \Omega(M)  \ar[r]^I &  
   W(\Hh)_{\com} \otimes \Omega(M) \ar[uu]^{ c \otimes \id }  
}
\end{equation}
There is a natural projection from $R $ onto $\gm $
given by $\pi (h_1,h_2,m)\equal (h_1^{-1} h_2, h_2^{-1}\cdot m)$. This endows $R \toto H \times M$
with a structure  of principal $H$-groupoid over $\gm \toto M$,
where $H$ acts on $R$ by  $h\cdot (h_1,h_2,m) \equal  (hh_1, hh_2, h \cdot m)$.

Therefore $\Omega(R_{\com})^{H-basic}\simeq \Omega(\gm_{\com})$.
Taking $H$-basic elements in the previous diagram, we obtain
\begin{equation}
\label{eq: tableau'}
\xymatrix{
\parallel{\Omega}(\gm_{\com})\parallel\ar[r]^{I}&{\Omega}(\gm_{\com})\\
{\big(}W({\Hh})\otimes\Omega(M)\big)^{H-basic}\ar[u]\ar[d]^{F^{H-basic}} &\\
{\big(}{\parallel}W({\Hh})_{\com}\parallel\otimes\Omega(M){\big)}^{H-basic}\ar[r]^I&
{\big(}W({\Hh})_{\com}\otimes\Omega(M){\big)}^{H-basic}\ar[uu]_{(c\otimes\id)^{H-basic}}\\
}
\end{equation}
Therefore we obtain two  equivalent 
descriptions of the natural chain map
$$K:   \big(W(\Hh)  \otimes \Omega(M)   \big)^{H-basic} \to \Omega(\gm_\com). $$ 

\begin{prop} 
\label{pro:5.17}
If $H$ is a  compact Lie group, then  $K$ is a quasi-isomorphism.
\end{prop}

To prove this result, we need some preliminaries.
Recall that a {\em $W(\Hh) $-module} $\Omega$ is an $H$-differential algebra endowed with an action $W(\Hh)\otimes \Omega  \to \Omega $
 which  is a homomorphism
of $H$-differential algebras. A {\em $W(\Hh) $-algebra} is a $H$-differential algebra
with a structure of  $W(\Hh) $-module. Note that an $H$-differential algebra
can be endowed with a structure of $W(\Hh)$-module if and only if it admits a connection.
  A $W(\Hh) $-module is said to be {\em acyclic}
if its cohomology vanishes in strictly positive degree and is ${\mathbb R}$
in degree $0$ (see \cite{guillemin}). 

According to  (Section 4.4 \cite{guillemin}), for any acyclic $W(\Hh)$-algebras $A$ and $B$,  
 the inclusion
$B \otimes \Omega \hookrightarrow A \otimes B \otimes \Omega$ 
induces an isomorphism 
\begin{equation}\label{eq:isoincl} H^*(\big( B \otimes \Omega\big)^{H-basic} ) \simeq 
H^*(\big( A \otimes  B \otimes \Omega\big)^{H-basic} ) . \end{equation} 
 Hence for any acyclic $W(\Hh)$-algebras $A$ and $B$,
an isomorphism $K(A,B)$ from $H^*(\big( A \otimes \Omega\big)^{H-basic} )$ 
to $H^*(\big( B \otimes \Omega \big)^{H-basic} )$  can be canonically constructed by
  \begin{equation}\label{eq:kab} K(A,B) \equal  H(j^{H-basic})^{-1} \smalcirc  H(\tau^{H-basic})
 \smalcirc H(i^{H-basic}) .\end{equation} Here,
\begin{itemize}
\item $\tau: A \otimes B\otimes 
\Omega \to B \otimes A\otimes \Omega $ is obtained by flipping $A$ and $B$
according to the Quillen rule, i.e.,
$$ \tau(a \otimes b \otimes \omega) = (-1)^{kl} b \otimes a \otimes \omega  \ \ \ \ \forall a \in A^k,b \in B^l, \omega \in \Omega $$
\item 
 $i :A \otimes \Omega \to B \otimes A \otimes \Omega $ is the canonical inclusion
$i(a \otimes \omega) \equal (1_B\otimes a \otimes \omega)$,
\item   $j :B \otimes \Omega \to A \otimes B \otimes \Omega $ is the canonical inclusion
 $j(b \otimes \omega) \equal (1_A\otimes b \otimes \omega)$.
\end{itemize}

\begin{lem}
 \label{lem:lemacyclicmodule}
 Let $\Omega$ be an $H$-differential algebra and $A, B$ 
acyclic $W(\Hh)$-algebras. Assume that $\phi: A\to B$
is a homomorphism of  $W(\Hh)$-algebras.
Then the homomorphism
 $$H(\big(\phi \otimes \id\big)^{H-basic}): H^*(\big( A \otimes \Omega\big)^{H-basic})
\to H^*(\big( B \otimes \Omega\big)^{H-basic})$$ is equal to $K(A,B)$
and is therefore an isomorphism. 
\end{lem}
\begin{pf}
 Let $\mbox{\j}_A$ be the homomorphism
of $G$-differential algebra $\mbox{\j}_A: W(\Hh) \longrightarrow A$  defined 
by $w \to w \cdot 1_A$.  As a first step, we show that 
$$H(\big( \mbox{\j}_A \otimes \id \big)^{H-basic}) : H^*(\big(W(\Hh) \otimes \Omega \big)^{H-basic} )\to H^*(\big(A \otimes \Omega \big)^{H-basic} )$$
 is an isomorphism.  

It is simple to check that the following diagram is commutative,
$$
\xymatrix{
 W(\Hh)\otimes \Omega \ar[r]^{\mbox{\j}_A \otimes \id} \ar[d]_{i_1}&
  A \otimes \Omega \ar[d]_{i_2}  \\
 A \otimes W(\Hh)\otimes \Omega    &  A \otimes A \otimes \Omega \\
      A    \otimes  \Omega  \ar[u]^{i_3} \ar[ru]_{i_4} &     \\
}
$$
where $i_1,i_2,i_3,i_4$ are defined by, $\forall w \in W(\Hh), a \in A ,\omega \in \Omega$
$$ \begin{array}{ccc} 
 i_1(w \otimes \omega)  &    \equal  &  (1_A \otimes w        \otimes \omega) \\
 i_2(a \otimes \omega)  &    \equal  &  (1_A \otimes a       \otimes \omega)\\
 i_3(a \otimes \omega)  &    \equal  &  (a   \otimes 1_{W(\Hh)} \otimes \omega) \\
 i_4(a \otimes \omega)  &    \equal  &  (a   \otimes 1_{A}       \otimes \omega) .  \\
 \end{array}$$
By Eq.~(\ref{eq:isoincl}), $i_1,i_2,i_3,i_4$, 
when being restricted to the $H$-basic elements,
 induce isomorphism of  cohomology.
Therefore $$H(\big(\mbox{\j}_A \otimes \id\big)^{H-basic}) \equal H(i_2^{H-basic})^{-1} 
\smalcirc H(i_4^{H-basic}) \smalcirc H(i_3^{H-basic})^{-1} \smalcirc H(i_1^{H-basic}) $$ is indeed an isomorphism.

The commutativity of the diagram
$$\comtri{W(\Hh)}{\mbox{\j}_A}{A}{\mbox{\j}_B}{\phi}{B}$$
implies the commutativity of the diagram
\begin{equation} \label{eq:commdernier}
\xymatrix{
  H^*(\big(W(\Hh) 
\otimes \Omega \big)^{H-basic})  \ar[rr]^{ H(\big( \mbox{\j}_A\otimes
 \id\big)^{H-basic}) }  \ar[rd]_{ H(\big( \mbox{\j}_B
   \otimes \id\big)^{H-basic}) } &   &   H^*(\big(A \otimes \Omega \big)^{H-basic}) 
\ar[ld]^{ H( \big( \phi \otimes 
\id \big)^{H-basic})  }  \\
&     H^*(\big(B \otimes \Omega \big)^{H-basic})& \\
}
\end{equation}
Since $H( \big( \mbox{\j}_A \otimes \id \big)^{H-basic})$ 
and  $H( \big( \mbox{\j}_B \otimes \id \big)^{H-basic})$   are isomorphisms, 
 the commutativity of diagram (\ref{eq:commdernier})
 implies that $H(\big(\phi \otimes \id \big)^{H-basic})$
 must be an isomorphism. Moreover,
if $\phi$ and $\psi $ are two such homomorphisms, then
\begin{equation}
\label{eq:iso}
H(\big(\phi \otimes \id \big)^{H-basic})
=H(\big(\psi \otimes \id \big)^{H-basic}).
\end{equation}

Now  consider  the following two homomorphisms
  of acyclic $W(\Hh)$-algebras from   $A$ to $A \otimes B$
given by  $i_1 (a)= a \otimes 1_B $
and $j_1 (a)= 1_A \otimes \phi(a)$, $\forall a\in A$.
Then  we have $i=i_1 \otimes \id$ and $\tau \smalcirc j \smalcirc
(\phi \otimes \id)=j_1\otimes \id$. By Eq. (\ref{eq:iso}),
we have
\be && H( i^{H-basic} )=H( \big(i_1 \otimes \id\big)^{H-basic})=
H( \big(j_1 \otimes \id\big)^{H-basic})\\
&&
=H( \tau^{H-basic} ) \smalcirc H( j^{H-basic} ) \smalcirc H( \big( \phi \otimes id\big)^{H-basic}). \ee
The conclusion thus follows from   Eq.~(\ref{eq:kab}) immediately.
\end{pf}

Now we can prove Proposition \ref{pro:5.17}.
\begin{pf} 
According to \cite{guillemin}, 
$\parallel \Omega((H\times H)_\com) \parallel$ (denoted
by  $\Omega({\mathcal E})$ in \cite{guillemin}) 
is a acyclic $W(\Hh)$-algebra. Lemma \ref{lem:lemacyclicmodule},
 applied to the case
that  $\Omega \equal \Omega(M)$, $A \equal W(\Hh )$ 
and  \\ $B \equal \parallel \Omega((H\times H)_\com) \parallel$,
implies that 
$$c_{\tilde{\zeta}}\otimes \id : W(\Hh) \otimes \Omega(M)\to  \parallel \Omega((H\times H)_\com)\otimes \Omega(M) \parallel$$ induces an isomorphism in cohomology when restricted to  basic elements.  

Now the cohomology of $\parallel \Omega((H\times H)_\com)\otimes \Omega(M) \parallel^{H-basic}$
is equal to the cohomology of $\Omega(\Gamma_\com) $.
 Proposition \ref{pro:5.17} then follows immediately from 
the fact that $I:\parallel \Omega(\gm_\com) \parallel \to \Omega(\Gamma_{\com})$
induces an isomorphism in cohomology according to Proposition \ref{prop:proprieteintegration}. \end{pf}

\begin{numrk} An explicit quasi-isomorphism between the Cartan model and  the
simplicial model of an equivariant cohomology  was also
constructed by Meinrenken  \cite{Mein}.
It would be interesting
to compare these two constructions.
\end{numrk}

Let $P \to M $ be an $H$-equivariant principal $G$-bundle.
In \cite{BT}, Bott-Tu
 introduced a Chern-Weil map with values in the Weil model
associated to any $H$-invariant connection $\theta$ on $P \to M$.
More precisely, they  constructed an
 $H$-basic connection on $ W(\Hh)  \otimes \Omega(P) $
by
 \begin{equation} \label{eq:connbotttu}
\Xi\equal  \sum_{i\equal 1 }^{dim(\Hh)} (1 \otimes \xi^i)
\otimes  L_i+1_{W(\Hh)} \otimes \theta
 \, \, \, \, \in W(\Hh) \otimes  \Omega(P) \otimes \Gg,
\end{equation}
where $L_i \equal  -  <\theta,\hat{X}_i>$ is a $\Gg$-valued function on $P$.
This connection, according to Proposition \ref{prop:propbasealgweil}, 
induces a map $ S(\Gg^* )^G \to 
 \big(   W(\Hh)   \otimes \Omega(M)   \big)^{H-basic} $
that we denote by $z_{BT}$.
On the other hand, $\theta \in \Omega(P)\otimes \Gg$ can be
considered as a pseudo-connection
of $ (W(\Hh)_{\com} \otimes \Omega(P))^{H-basic}$. 
The construction of Section \ref{sec:CWgroupoides} 
induces a map $z_{\theta}:S^*(\Gg)^G \to Z^* (\gm_\com)$.

\begin{them}
\label{theo:botttu}
The following diagram  commutes
  $$\comtri{S(\Gg^* )^G}{z_{BT}}{Z^*(\big( W(\Hh)  \otimes \Omega(M) \big)^{H-basic})}{z_\theta}{K }{Z^*(\gm_\com)  }  $$
\end{them}

We will need the following lemma first.

\begin{lem} \label{lem:verifier}
Assume that $\theta \in \Omega^1(P) \otimes \Gg$
 is an
 $H$-invariant connection on the principal $G$-bundle $P \to M$.
Let $\pi: H \times P \to P$ be the map defined by
$\pi(h,p)\equal  h^{-1}\cdot p $. 
We have 
$$ \pi^* \theta \equal   \sum_{i\equal 1 }^{dim(\Hh)}
\zeta^i  \otimes L_i +  1  \otimes \theta \,\,\,\,\,\, \in \Omega(H)
\otimes \Omega(P)\otimes \Gg \subset  \Omega(H \times P)\otimes \Gg .$$
\end{lem}
\begin{pf}
The proof is  a direct verification and is left to the reader.
\end{pf}

Now we are ready to  prove Theorem \ref{theo:botttu}.

\begin{pf} The $1$-form $\Xi$, as defined in Eq.~(\ref{eq:connbotttu}),
can be considered
as a pseudo-connection
on the $G$-differential simplicial algebra
$   W(\Hh)_{\com}   \otimes \Omega(P)  $.  
By definition of  $\psi_\eta$, we have 
$$ (\psi_\eta \otimes \id) (\Xi)
\equal  \sum_{i=1}^{dim({\mathfrak h})} \psi_\eta(1 \otimes \xi_i)  \otimes L_i +
1 \otimes \theta  \equal    \sum_{i=1}^{dim({\mathfrak h})}   \zeta_i  \otimes L_i +
1  \otimes \theta   .
$$
By Lemma \ref{lem:verifier},  we have
$ (\psi_\eta \otimes \id)(\Xi ) \equal   \pi^* \theta \in \Omega( R_\com)$,
where $R $ denotes  the direct product of the pair groupoid
$H \times H \toto H$
with the manifold $P \toto P$.

 By the functoriality properties of the Chern-Weil map, we obtain the 
commutative diagram
 \begin{equation} \label{eq:cont2} \begin{array}{ccccc} 
  &  &  \parallel \Omega(R_\com)  \parallel & \stackrel{I}{\to}
& \Omega(R_\com)\\
    & c_{\tilde{\pi^* \theta}}\nearrow & &  & \\
W(\Gg) & & & & \rdiagup{ \psi_\eta \otimes \id}  \\
  &c_{\tilde{\Xi}}\searrow &  & &  \\
   & &  \parallel W(\Hh)_\com \parallel \otimes \Omega(P)  & \stackrel{I}{\to}  &        W(\Hh)_\com  \otimes \Omega(P)  \\
\end{array}  \end{equation}
Since the image of $\Xi$ under $F$ 
is $ \tilde{\Xi}$, which is a connection
for the $G$-differential algebra
 $ \parallel W(\Hh)_\com \parallel \otimes \Omega(P)$,
 we also have a commutative diagram
\begin{equation}\label{eq:cont3} \comtri{W(\Gg)}{c_{\Xi}}{W(\Hh)\otimes \Omega(P) }{c_{\tilde{\Xi}}}{T}{ \parallel W(\Hh)_\com \parallel \otimes \Omega(P)}   \end{equation}
Combining Eq.~(\ref{eq:cont2}) and Eq.~(\ref{eq:cont3}), we obtain
 the commutative diagram
\begin{equation} \label{eq:cont4} 
\xymatrix{ 
   &\parallel \Omega(R_\com)  \parallel \ar[r]^I &  \Omega(R_\com)  \\
 W(\Gg)  \ar[ru]^{  c_{\tilde{\pi^* \theta}} } \ar[r]^{c_{\Xi }} \ar[rd]^{ c_{\tilde{\Xi}}  }  &
 W(\Hh) \otimes \Omega(P) \ar[d]^{\Xi} &   \\
    &  \parallel W(\Hh)_\com \parallel \otimes \Omega(P)  \ar[r]^I & W(\Hh)_\com  \otimes \Omega(P) 
\ar[uu]^{ \psi_{\eta} \otimes \id }  \\
}
\end{equation}
Since all the arrows of this diagram are  $G$-$H$-bimodule maps
($W(\Gg)$ is considered as a trivial $H$-module), we can restrict ourself to elements
which are both $G$- and $H$-basic.
We then obtain the following commutative diagram:

\begin{equation} \label{eq:cont5} 
\xymatrix{ 
   &\parallel \Omega(\Gamma_\com)  \parallel \ar[r]^I &  \Omega(\Gamma_\com)  \\
 S(\Gg^*)^G \ar[ru]^{  c_{\tilde{ \theta}} } \ar[r]^{z_{BT}} \ar[rd]^{ c_{\tilde{\Xi}}  }  &
 {\big(} W(\Hh) \otimes \Omega(P){\big)}^{H-basic} \ar[d]^{T^{H-basic}} &   \\
    & {\big(} \parallel W(\Hh)_\com \parallel \otimes \Omega(P) {\big)}^{H-basic} \ar[r]^I & 
{\big(} W(\Hh)_\com  \otimes \Omega(P) {\big)}^{H-basic}
\ar[uu]^{ (\psi_{\eta} \otimes \id)^{H-basic} }  \\
}
\end{equation}

Now the composition of $c_{\tilde{\theta}}$  and $I$ is  $z_{\theta}$,
the  Chern-Weil map on the cochain level,
while  $K$ is  the composition  $(\psi_{\eta} \otimes \id)^{H-basic} 
\smalcirc I \smalcirc T^{H-basic}$.  Therefore we have 
$ z_{\theta} \equal  K \smalcirc z_{BT}  .$\end{pf}

\begin{rmk}
For equivariant Chern-Weil map in the Cartan model, we 
refer the reader to \cite{BGV, Jeffrey}.
\end{rmk}

\end{document}